\documentclass[12pt]{article}
\usepackage{amssymb}
\usepackage{amsfonts}
\usepackage{amsmath}
\usepackage{hyperref}
\usepackage{xcolor}
\usepackage{yfonts}
\usepackage{graphicx,epstopdf}
\usepackage{tikz}

\usepackage{tikz-3dplot}

\usepackage{lscape}

\usepackage{xcolor}
\usepackage{yfonts}
\usepackage{hyperref}
\usepackage{float}

\usetikzlibrary{shapes.geometric, calc}
\usetikzlibrary{decorations.markings,shapes.arrows, shapes.symbols}
\usetikzlibrary{decorations.pathreplacing}

\newtheorem{theorem}{Theorem}[section]
\newtheorem{proposition}[theorem]{Proposition}
\newtheorem{lemma}[theorem]{Lemma}
\newtheorem{corollary}[theorem]{Corollary}

\newtheorem{definition}[theorem]{Definition}
\newtheorem{example}[theorem]{Example}
\newtheorem{remark}[theorem]{Remark}

\newenvironment{proof}[1][Proof]{\textbf{#1.} }{\ \rule{0.5em}{0.5em}}

\begin{document}

\author{ Jingyan Li\\
Beijing Institute of\\ Mathematical 
Sciences and \\Application,
Beijing, China\\
\and 
Yuri Muranov \\
University  of  Warmia and \\
Mazury in Olsztyn, Olsztyn\\
Poland \\
\and Jie Wu\\
Beijing \\ 
Institute of Mathematical \\
Sciences and Application,\\
Beijing, China\\
\and Shing-Tung Yau\\
Tsinghua University, \\
Beijing, China\\
}
\title{Primitive path  homology }
\date{November, 2024}
\maketitle

\begin{abstract} In this paper we introduce a primitive path homology theory on the category of simple digraphs. On the subcategory of  asymmetric digraphs, this theory coincides  with the path homology theory which was introduced by Grigor'yan, Lin, Muranov,  and Yau, but these theories are different  in general case.   We study properties of the primitive path homology and  describe relations between the primitive path  homology and the path homology. Let  $a,b$  two different vertices of a digraph.   Our  approach gives a possibility to  construct  primitive    homology theories  of  paths which have a given tail vertex $a$ or (and) a given  head vertex $b$.  We study these theories and describe also relationships between them  and the  path homology theory. 
\end{abstract}

\bigskip

\tableofcontents

\section{Introduction}\label{S1}
\setcounter{equation}{0}

In this paper, we develop the path homology theory  (GLMY theory)  that was created by  Grigor'yan, Lin, Muranov, and Yau in  \cite{Path}, \cite{Mi2012}, \cite{MiHomotopy}, \cite{Mi2}.   The path homology theory closely relates to other homology theories on discrete objects. The simplicial homology of a simplicial complex is isomorphic to the path homology of a digraph which is functorially defined by the   complex   \cite{Mi3}. The deep  relations between the path cohomology theory of digraphs, and Hochschild cohomology of algebras were obtained in \cite{MiHH}. The category of discrete spaces coincides with the category of transitive digraphs, and hence, Alexandroff homology  \cite{Pal}  gives  a homology  of transitive digraphs which coincides with the path homology  \cite{Muranov}.  The singular cubic homology coincides with the path homology in the category of transitive digraphs and in the category of cubical digraphs   \cite{MurJim23}, \cite{Mi5} but these  homology theories   are different   in general \cite[Pr. 13, Pr. 14]{Mi4}. 

Recently, a cell approach was developed to study path homology in \cite{XinFu}, \cite{Sasha}, \cite{Yau15}, \cite{Yau23}, \cite{Yau24}, and the path homology of cycles, circuits and closed walk was constructed in \cite{MurJC} and \cite{Circuits}.
In both approaches, paths with some restrictions on the tail and head vertex arise naturally.

 In this paper,   we introduce a primitive path homology theory on the category of simple digraphs. This theory coincides with the path homology theory on the subcategory of simple asymmetric digraphs and differs from the path homology theory in general.  We study properties of the primitive path homology and  describe relations between this homology and  the path homology. Let  $a,b$ be  two different vertices of a digraph.   Our approach gives a possibility to construct primitive homology theories of paths that have a given tail vertex $a$ or (and) a given head vertex $b$.  We also describe relations between these theories,  the path homology and the path homology of clusters which were defined in \cite{Sasha}. We give a number of examples that illustrate the results obtained. The paper is organized as follows.

In Section \ref{S2},  we recall the basic definitions of the digraph theory
and introduce the path homology theory of  digraphs.

In Section \ref{S3}, we define  and study  the primitive path homology theory of digraphs. We describe relations between  primitive path homology   and the path homology and provide several examples. In particular, we prove that the primitive path homology coincides with the path homology on the category of asymmetric digraphs  and give  an example which  shows that these theories differ in general.

In Section \ref{S4},  we  define a cluster subgraph $G^{[a,b]}$ of a digraph $G$.  The cluster digraph  $G^{[a,b]}$  in $G$ is a  subgraph given   by all vertices and all  arrows  of paths with  a tail vertex $a$ and a head vertex $b$ in the digraph $G$. We describe relations between primitive homology of cluster  digraphs and primitive homology of ambient digraphs. Similarly, we define tail fixed digraphs, and head fixed digraphs which are subgraphs of a digraph $G$  and describe its relations to  the primitive path homology of cluster digraphs and ambient digraphs.

In Section \ref{S5},  for a digraph $G=(V,E)$ and two vertices $a,b \in V$,   we define a primitive homology of paths which have the tail vertex $a$  and  the  head vertex $b$. We call 
this homology by \emph{primitive  $(a,b)$-cluster homology} of the digraph $G$.  We describe basic properties of this homology  and 
 its relation to primitive homology of the cluster digraph in $G$  and to  primitive homology of $G$.  We  provide also several examples.  

In Section \ref{S6}, we  introduce the  primitive homology of paths with a fixed tail vertex  and the  primitive homology of paths with a fixed  head vertex. We  describe relations between them and its relations to the various primitive homology which are considered in Section \ref{S5}. We  provide also several examples.

In Section \ref{S7}, we prove functoriality of primitive  homology  theories constructed in Sections  \ref{S5} and \ref{S6} on  the category $\mathcal{AD}$ of asymmetric digraphs. 

\section*{Acknowledgments}
This work is supported by the start-up research funds of the first author and the third author from BIMSA. The third author was supported by High-Level Scientific Research Foundation of Hebei Province.

\section{The path homology}\label{S2}
\setcounter{equation}{0}

In this section, we recall  the necessary  notions of the graph theory  and  the path homology theory  \cite{Gary}, \cite{Mi2012},  \cite{MiHomotopy}, \cite{Mi3}, \cite{MiHH},  \cite{Sasha}.

In this paper, we only consider  \emph{simple digraphs} $G=(V_G,E_G)$ with a nonempty and finite set $V_G$ of \emph{vertices} and with a set $E_G$ of 
ordered pairs $(v,w)\in E_G$ of different vertices $v,w\in V$. A pair $(v,w)\in E_G$  is  called an \emph{arrow}, denoted $v\to w$, where $v$ is the \emph{tail} and $w$ is the \emph{head} of the arrow. In this case we write $v=t(v\to w)$ and $w=h(v\to w)$.  For two vertices $v,w\in V_G$ the arrows $v\to w, w\to v\in E_G$ are 
called \emph{symmetric}. A  digraph $G$ is called \emph{asymmetric} 
if it doesn't contain symmetric arrows.

A digraph $G=(V_G,E_G)$ is a \emph{subgraph} of a digraph 
$H=(V_H,E_H)$,  and we write $G\sqsubset H$,   if $V_G\subset V_H$ and $E_G\subset V_H$. A subgraph  $G\sqsubset H$ is called the \emph{induced subgraph} and is denoted $G\subset H$, if for any $v, w\in E_G\subset E_H$ and $(v\to w)\in E_H$  implies $(v\to w)\in E_G$. Any subgraph $G\sqsubset H$ uniquely defines an  minimal induced  subgraph $\overline{G}\subset H$ containing $G$. 

Let $G$ and $H$ be two digraphs. A \emph{map}  $f\colon G\to H$ is 
given  by a map  $f\colon V_{G}\rightarrow V_{H}$ of vertices  such that for every arrow $(v\to w)\in E_G$ we have either $(f(v)\to f(w))\in E_H$ or $f(v) =f(w)\in V_H$.  A map $f\colon G\to H$ is  a  \emph{homomorphism} if 
$(f(v)\to f(w))\in E_H$ for any $(v\to w)\in E_G$.

Thus we obtain the category $\mathcal D$ of digraphs and digraph maps. By a natural way,  there are  two subcategory  of  $\mathcal D$ defined.  We have a  subcategory $\mathcal{ND}\subset \mathcal D$ with the same set of objects as in $\mathcal D$  and with the morphisms which are homomorphisms and 
 a  subcategory $\mathcal{AD}\subset \mathcal D$ with the same set of morphisms  as in $\mathcal D$  and with the set of objects consisting of asymmetric digraphs.  

Recall the definition of path homology groups of a digraph $G$ with coefficients in a unitary commutative ring $R$. Let $V$ be a finite set whose elements we call \emph{vertices}. For $n\geq 0$, an \emph{elementary} $n$-\emph{path} in $V$  is any sequence $i_{0},\dots,i_{n}$ of vertices.   Let  
$\Lambda_n=\Lambda _{n}( V, R)$ be a free $R$-module generated by the elements $e_{i_{0}\dots i_{n}}$ where $i_{0},\dots,i_{n}$ is an elementary $n$-path in $V$.  
The elements of $\Lambda _{n}$ are called $n$-\emph{paths} in  $V$. We set also $ \Lambda _{-1}=0$. For $n\geq 1$, define the homomorphism  (\emph{differential}) $\partial\colon \Lambda _{n}\rightarrow \Lambda _{n-1}$ by setting on   basic elements 
\begin{equation}\label{2.1}
 \partial (e_{i_{0}...i_{n}})=\begin{cases} \sum\limits_{m=0}^{n}\left( -1\right)^{m}e_{i_{0}\dots \widehat{i_{m}}\dots i_{n}}& \text{for} \ n\geq 1,\\
 0& \text{for} \  n=0.
\end{cases}
\end{equation}
where $\widehat{i_{m}}$ means omitting the vertex $i_m$. For $n\geq 0$, we can write down the differential  in the form
$$
 \partial (e_{i_{0}...i_{n}})=\sum\limits_{m=0}^n(-1)^m \partial^m(e_{i_{0}...i_{n}}) 
$$
where $\partial^m\colon \Lambda _{n}\rightarrow \Lambda _{n-1}$ is a homomorphism given on the basic elements by 
\begin{equation}\label{2.2}
\partial^m(e_{i_{0}...i_{n}})=\begin{cases} e_{i_{0}\dots \widehat{i_{m}}\dots i_{n}} & \text{for} \  n\geq 1,  \,  0\leq m\leq n, \\
0& \text{for} \ m=n=0.\\
\end{cases}
\end{equation}
The differential $\partial$  has the property $\partial\circ \partial=0$  and the modules $\Lambda_n$  with the differential $\partial$ form a chain complex $\Lambda_*=\Lambda_*(V, R)$.

 We call an elementary path $e_{i_{0}\dots i_{n}}$ \emph{irregular} if the path 
$i_{0}, \dots,  i_{n}$
  has two equal consecutive vertices  and \emph{regular} otherwise. 

For $n\geq 1$, let $I_{n}=I_n(V,R)$ be a submodule of $\Lambda _{n}( V, R)$ that is generated by all irregular elementary paths. Set $I_0=I_{-1}=0$. For $n\geq -1$,  denote by 
$\mathcal R_n=\mathcal {R}_{n}(V, R)$  the  quotient module 
$\Lambda _{n}/I_{n}=\Lambda_n(V, R)/I_n(V, R)$. The elements of modules $\mathcal R_n$ we will continue to call  \emph{regular paths in the set $V$}.

Since  $\partial (I_{n})\subset I_{n-1}$, the modules  $\mathcal R_n$ with the induced  differential  $\partial$, which we continue to denote $\partial$, form a chain  complex 
$\mathcal R_*$.

Let $G=(V_G, E_G)$ be a digraph. For $n\geq 1$, an elementary path $e_{i_{0}\dots i_{n}}$ in $V_G$ is called \emph{allowed} if $(i_{k}\rightarrow i_{k+1})\in E_G$ for $0\leq k\leq n-1$. We note that  every allowed path is regular since we consider only simple digraphs.  The \emph{length} of an allowed  path $e_{i_{0}\dots i_{n}}$ equals $n$ and it is denoted by  $|e_{i_{0}\dots i_{n}}|$. 
Every path $e_{i_0}$ is allowed and has the length $0$.   For an allowed path $e_{i_{0}...i_{n}}$ we call the vertex $i_0=t(e_{i_{0}...i_{n}})$ by \emph{the tail of the path} and the vertex 
$i_n=h(e_{i_{0}...i_{n}})$ by \emph{the head of the path}. The \emph{oriented distance} 
$d(v,w)$ between two vertices $v,w\in E_G$ is the minimal length of an allowed path with  the tail $v$ and the head $w$. If such a path doesn't exist then the  distance equals $\infty$.

For $n\geq 0$, let $\mathcal A_{n}=\mathcal A_{n}(G,R)$ be a module that is generated by all allowed elementary $n$-paths and let  $\mathcal A_{-1}=0$. We  set also $\mathcal A_n=0 \,  (n\geq 1)$  if elementary allowed $n$-paths  do not exist in $G$. 
We will continue to call elements of $\mathcal A_n(G)$ by \emph{allowed paths}. 
 This  module can be naturally considered as  a submodule of  
$\mathcal R_{n}(V_G,R)=\mathcal R_n(V_G)$. The inclusion $\eta\colon \mathcal A_n(G)\to \mathcal R_n(V_G)$ is given  on the allowed elementary paths $e_{i_0\dots i_n}$, which form the basis of $\mathcal A_n(G)$,   by $\eta(e_{i_0\dots i_n})=[e_{i_0\dots i_n}]$ where the classes 
$[e_{i_0\dots i_n}]$ give a basis of  $\eta(A_n(G))\subset \mathcal R_n(G)$. 
For $n\geq 0$, define a submodule $\Omega_n(G,R)\subset \mathcal R_{n}(V_G,R)$ by setting 
\begin{equation}\label{2.3}
\Omega_n=\Omega _{n}(G)=\Omega _{n}( G, R)=\left\{ w\in {\mathcal A}
_{n}(G)\, |\, \partial w\in {\mathcal A}_{n-1}(G)\right\},
\end{equation}
where $\partial$ is the differential 
$\partial \colon \mathcal R_n(V_G)\to \mathcal R_{n-1}(V_G)$. 
We set $\Omega_{-1}=0$. It follows from (\ref{2.3}) that $\partial(\Omega_n)\subset \Omega_{n-1}$ for $n\geq 0$. Hence, the modules $\Omega_n$ with the differential $\partial$ form  a  chain complex  $\Omega_*$ which we
 call a \emph{path chain complex} of the digraph $G$.  

The  homology groups $H_*\left(\Omega_*\right)$ are called \emph{path homology groups of the digraph $G$} and are denoted by $H_*(G)=H_*(G, R)$. 

A digraph map $f\colon G\rightarrow H $ induces  a morphism of chain complexes $f_*\colon \Omega_*(G,R)\to \Omega_*(H,R)$ which on the set of allowed elementary paths 
is given by the condition 
$$
f_*(e_{i_0\dots i_n}) =\begin{cases} e_{f(i_0)\dots f(i_n)} & \text{ if } \  f(i_j)\ne f(i_{j+1}) \ \text{for any} \ 0\leq  i_j\leq n-1, \\
        0& \text{otherwise}.\\
\end{cases}
$$
We obtain an induced homomorphism  $f_*\colon H_*(G,R)\to H_*(H,R)$.
Thus, the path homology  is functorial on  the categories 
$\mathcal D$,  $\mathcal {ND}$, and $\mathcal {AD}$.

Recall  the notion of an augmentation  of a  chain complex 
$C_*=(C_n, \partial)  \, (n\geq 0)$ of $R$-modules with a differential $\partial \colon C_n\to C_{n-1}$. An \emph{augmentation} of $C_*$ is given by a homomorphism 
\begin{equation}\label{2.4}
\varepsilon\colon C_0\longrightarrow  R \ \ \text{such that} \ \ \varepsilon \partial=0. 
\end{equation}
Thus we obtain an  \emph{augmented chain complex} 
\begin{equation}\label{2.5}
\dots \overset{\partial}\longrightarrow C_2\overset{\partial}\longrightarrow C_1
\overset{\partial}\longrightarrow C_0 \overset{\varepsilon}\longrightarrow R\longrightarrow 0. 
\end{equation} 
which we denote $\widetilde C_*$.
\begin{definition}\label{d2.1} \rm 
 The \emph{reduced homology}  $\widetilde{H} (C_*)$ of a chain complex  $C_*$   is the homology of the augmented chain complex   $\widetilde C_*$ in (\ref{2.5}).
\end{definition}

For every set $V$, the chain complex $\Lambda_*(V, R)$ has the natural augmentation 
\begin{equation}\label{2.6}
\varepsilon\colon \Lambda_0(G,R)\longrightarrow  R,  \ \ 
 \varepsilon \left(\sum_i r^ie_i \right) = \sum_i r^i \ \ \text{where} \ \ r^i\in R, \ i\in V.  
\end{equation} 
For a digraph $G$, the augmentation of the chain complexes 
$\mathcal R_*(V_G,R)$ and $\Omega_*(G, R)$ is defined similarly (\ref{2.6}). 
Hence, for every digraph $G$,  the \emph{augmented path homology groups} $\widetilde{H}_*(G)=\widetilde{H}_*(G,R)$  are defined.

Let  $G$ be  a digraph and  $w\in\mathcal A_n(G, R)$ be an allowed  path. For $n\geq 1$,  we can write 
$w$ in the form  
\begin{equation}\label{2.7}
w=\sum_{a,b\in V}w^{[a,b]}, \ \ \  
w^{[a,b]} =\sum_{i_0=a, i_n=b}r^{i_0i_1\dots i_{n-1}i_n}e_{i_0i_1\dots i_{n-1}i_n}
\end{equation}
where   $r^{i_0i_1\dots i_{n-1}i_n}\in R$ and $e_{i_0i_1\dots i_{n-1}i_n}$  is an elementary allowed  path with 
$$
a=i_0=t\left(e_{i_0i_1\dots i_{n-1}i_n}\right),  \ \  b=i_n=h\left(e_{i_0i_1\dots i_{n-1}i_n}\right).
$$
We set $w^{[a,b]}=0$ if there are no elementary paths fitting into $w$ with the tail vertex $a$ and the head vertex $b$. 
For $n=0$  and  $a\ne b$ we have $w^{[a,b]}=0$. For   $a=b$ we have  $w^{[a,a]}=0$ for $n=1$ and $w^{[a,a]}=r^ae_a\, (a\in V)$ for $n=0$. Hence, for $a=b$   we can write every path $w\in \mathcal A_0(G)$ in the form similar  to (\ref{2.7})
\begin{equation}\label{2.8}
w=\sum_{a\in V}w^{[a,a]}=\sum_{a\in V}r^{a}e_{a}, \ \ r^{a}\in R. 
\end{equation}

\begin{definition}\label{d2.2} \rm  For $n\geq 1$, a  path $w^{[a,b]}\in \mathcal A_{n}(G)$  in (\ref{2.7}) and for $n=0$ a path  $w^{[a,a]}\in \mathcal A_{0}(G)$ in  (\ref{2.8})  is called a \emph{cluster path in dimension} $n\geq 0$
(see also  \cite[\S 2.2]{Sasha}). The cluster path  is \emph{elementary} if it is given by an elementary path. 
\end{definition} 

 Now we slightly reformulate the result from   \cite[Lemma 2.2]{Sasha}. 
\begin{theorem}\label{t2.3} Let  $G$ be a digraph and $w \in \mathcal A_n(G)\, (n\geq 0)$.  Then $w\in \Omega_n(G)$ if and only  if $w^{[a,b]}\in \Omega_n(G)$ for every   $w^{[a,b]}$ fitting into the decomposition $w$ in  (\ref{2.7}) or in (\ref{2.8}). 
\end{theorem}
\begin{proof}  For $n=0$ it is nothing to prove. For $n\geq 1$,   if  $w^{[a,b]}\in \Omega_n(G)$ then
$\partial(w^{[a,b]})\in \mathcal A_{n-1}(G)$ and 
$$
\partial(w)=  \partial\left(\sum_{a,b\in V}w^{[a,b]}\right)=\sum_{a,b\in V}\partial(w^{[a,b]})\in \mathcal A_{n-1}. 
$$ 
Hence $w\in \Omega_n(G)$. The inverse statement is proved in  \cite[Lemma 2.2]{Sasha} 
\end{proof} 

By technical reasons, we set  $\Omega_n^{[a,b]}(G)=0$ if there are no $(a,b)$-cluster paths   in dimension $n$.

\begin{corollary}\label{c2.4} For every digraph $G=(V,E)$ and $n\geq 0$, the module $\Omega_n(G)$  is a direct sum of submodules 
$$
\Omega_n(G)=\bigoplus_{a,b\in V }  \Omega_n^{[a,b]}(G) 
$$
where
$$
\Omega_n^{[a,b]}(G)=\langle w^{[a,b]}\in \Omega_n(G)\, | \, a,b\in V\rangle
$$
is  generated by 
$(a,b)$-cluster paths  in dimension $n$.   
\end{corollary}

\section{Primitive path homology}\label{S3}
\setcounter{equation}{0}

In this Section we define  a primitive path homology   of digraphs and describe its relation to the path homology. We prove that  these theories  coincide on the category $\mathcal{AD}$ of asymmetric digraphs.  We describe algebraic properties of the primitive homology and provide an example which illustrate the difference between the primitive path homology and the path homology.

Let $G$ be a digraph and $R$ be a ring  of coefficients.  We note that for $n\geq -1$  the module 
$\mathcal A_n(G)$ is a submodule of 
 $\Lambda_n(V_G)$. For $n\geq 0$, define a submodule 
\begin{equation}\label{3.1}
\Pi_n=\Pi_n(G)=\left\{ w\in {\mathcal A}
_{n}(G)\, |\, \partial w\in {\mathcal A}_{n-1}(G)\right\}
\end{equation}
where $\partial \colon \Lambda_n(V_G)\to \Lambda_{n-1}(V_G)$ is the differential in (\ref{2.1}). We set $\Pi_{-1}(G)=0$. 
It follows immediately from (\ref{3.1}) that $\partial (\Pi_n)\subset \Pi_{n-1}$ for $n\geq 0$ and, hence,  the modules $\Pi_n(G)$ with the restriction of the differential $\partial$ form  a  chain complex $\Pi_*=\Pi_*(G)$. By  (\ref{3.1}) and (\ref{2.3},   for any digraph $G=(V,E)$ we have 
\begin{equation}\label{3.2}
\Pi_0(G)=\Omega_0(G)=\langle e_i\, | \, i\in V\rangle, \ \  \Pi_1(G)=\Omega_1(G)=\langle e_{ij}\, | \, (i\to j)\in E\rangle. 
\end{equation} 

\begin{definition}\label{d3.1} \rm The homology groups of the chain complex 
$\Pi_*(G)$ are called \emph{primitive path homology groups} with the coefficients $R$ of a digraph $G$ and are denoted by $\mathcal H_n(G)=\mathcal H_n(G,R)\, (n\geq 0)$.
\end{definition} 

\begin{example}\label{e3.2} \rm   Consider the complete digraph $G=(V,E)$ with the set of vertices    $V=\{0,1\}$ and  the set of  arrows $E=\{0\to 1, 1\to 0\}$. Let $R=\mathbb Z$. The digraph $G$ is homotopy equivalent to the one-vertex digraph  and,  hence, 
$$
H_n(G)=\begin{cases} \mathbb Z& \text{for} \ \ n=0, \\
                                      0&  \text{otherwise} 
\end{cases} 
$$
since the path homology  is homotopy invariant \cite{MiHomotopy}. 
Now we compute  primitive path  homology groups of $G$.  By (\ref{3.2}) we have $\Pi_0(G)=\mathcal A_0(G) =\langle e_0, e_1\rangle$ and $\Pi_1(G)=\mathcal A_1(G) =\langle e_{01}, e_{10}\rangle$.   For $n\geq 2$,  we have 
\begin{equation}\label{3.3}
\mathcal A_{n}(G)=\begin{cases} \langle e_{\tiny  \underbrace{01 \dots 10}_{n+1}},\, e_{\tiny\underbrace{10 \dots 01}_{n+1}}\rangle & \text{for} \  n=2k, \\
\langle e_{\tiny \underbrace{01 \dots 01}_{n+1}},\, e_{\tiny\underbrace{10 \dots 10}_{n+1}}\rangle & \text{for}  \  n=2k+1.\\
\end{cases}
\end{equation} 
Now we define the groups  $\Pi_n(G)$ for $n\geq 2$. For $n=2k\geq 2$,  the differential $\partial \colon \Lambda_{n}(G)\to \Lambda_{n-1}(G)$ on  the  basic element
$e_{\tiny \underbrace{01 \dots 10}_{n+1}}\in \mathcal A_n(G)$ contains the non-canceled summand 
$\partial^1(e_{\tiny \underbrace{01 \dots 10}_{n+1}})=(-1) e_{\tiny \underbrace{001 \dots 10 }_{n}}\notin \mathcal A_{n-1}(G)$.  Similarly,  for $n=2k+1$  the differential on the  basic element
$e_{\tiny\underbrace{10 \dots 01}_{n+1}}$ contains the non-canceled summand 
$\partial^1(e_{\tiny \underbrace{10 \dots 01}_{n+1}})=(-1) e_{\tiny\underbrace{001 \dots 01}_{n}}\notin \mathcal A_{n-1}(G)$. Hence, the image of the differential on every linear combination of the basic elements of $\mathcal A_{n}(G)$ doesn't belong to 
$\mathcal A_{n-1}(G)$ for $n\geq 2$  and   $\Pi_{n}(G)=0$ in these cases. 
The differential $\partial \colon \Pi_0(G)\to \Pi_{-1}(G)=0$ is trivial 
and the differential  $\partial \colon \Pi_1(G)\to \Pi_{0}(G)$ 
is given on the basic elements by equalities 
$\partial(e_{01})=e_1-e_0, \ \  \partial(e_{10})=e_0-e_1$. Hence 
$$
\mathcal H_n(G)=\begin{cases} \mathbb Z& \text{for} \ \ n=0, 1,  \\
                                      0&  \text{otherwise.} 
\end{cases} 
$$
Thus, the groups $H_n(G)$ and $\mathcal H_n(G)$ are not isomorphic in general case. 
\end{example}

Fix a ring of coefficients $R$.
Consider  a digraph map $f\colon G\rightarrow H $. Recall that  for any arrow $(v\to w)\in E_G$ we have $f(v)=f(w)\in V_H$ or $(f(v)\to f(w))\in E_H$.  For $n\geq 0$, define   a homomorphism of  modules  $$
f_{\sharp}\colon \Lambda_n(V_G)\to \Lambda_n(V_H)
$$
 setting on basic elementary paths  
\begin{equation}  \label{3.4}
f_{\sharp}\left( e_{i_{0}\dots i_{n}}\right) =
\begin{cases}
e_{f(i_{0})\dots f(i_{n})}, & \text{if }\ e_{f(i_{0})\dots f(i_{n})}\text{ is
regular,} \\ 
0, & \text{otherwise.}
\end{cases}
\end{equation}
We set also  $f_{\sharp}=0$ for $n=-1$. We note that the restriction of $f_{\sharp}$ defines a homomorphism 
$\mathcal A_{n}\left(G\right) \rightarrow \mathcal A_{n}(H)$ which we continue to denote $f_{\sharp}$. 

\begin{lemma}\label{l3.3} Let $f\colon G\to H$ be a digraph map and digraph $H$ is asymmetric. Then for every $n\ge 0$,  the following diagram is commutative
\begin{equation}\label{3.5}
\begin{matrix}
\mathcal A_{n}\left(G\right) &\overset{f_{\sharp}}\rightarrow &\mathcal A_{n}\left(H\right)\\
\\
\partial \downarrow && \downarrow \partial \\
\\
\Lambda_{n-1}\left( V_G\right) &\overset{f_{\sharp}}\rightarrow &\Lambda_{n-1}\left( V_H\right).  \\
\end{matrix} 
\end{equation} 
\end{lemma} 
\begin{proof} For $n=0$ it is nothing to prove since both maps $\partial$ are trivial.  Let $n\geq 1$ and 
$e_{i_0\dots i_n}\in \mathcal A_n(G)$.  We consider the following  cases:

(i)  the path $e_{f(i_0)\dots f(i_n)}$  is  not allowed, 

 (ii) the  path $e_{f(i_0)\dots f(i_n)}$  is  allowed. 

Consider the case (i). Let the path $e_{f(i_{0})\dots f(i_{n})}$ is not allowed. Then there is at least 
one pair of   consequent vertices $i_k,i_{k+1}$ in the path $e_{i_0\dots i_n}$ such that  $f(i_k)= f(i_{k+1})\in V_H$. Thus, the path  $e_{f(i_{0})\dots f(i_{n})}$ is 
not regular  and,  by (\ref{3.4}), 
\begin{equation}\label{3.6}
\partial f_{\sharp}(e_{i_0\dots i_n})=\partial (0)=0.
\end{equation}
 Furthermore, we have 
\begin{equation}\label{3.7}
\begin{split}
f_{{\sharp}}\partial(e_{i_0\dots i_n})&=f_{\sharp}\left(\sum_{m=0}^n (-1)^m e_{i_0\dots \widehat{i_m}\dots i_n}\right) \\
& = f_{{\sharp}}\left(\sum_{\underset{m\ne k, k+1}{m=0}}^n (-1)^m e_{i_0\dots \widehat{i_m}\dots i_n}\right)\\
&+  (-1)^k e_{f(i_0)\dots \widehat{f(i_k)}\dots f(i_n)}  \\ &+(-1)^{k+1} e_{f(i_0)\dots \widehat{f(i_{k+1})}\dots f(i_n)}.
\\
\end{split} 
\end{equation}
Every elementary path in the sum of the second row in (\ref{3.7}) has two  consequent vertices $i_k, i_{k+1}$ such that  $f(i_k)= f(i_{k+1})\in V_H$. Hence,  by the definition $f_{\sharp}$,  the second row in (\ref{3.7}) is zero. The paths of the two bottom rows in (\ref{3.7}) differ only by sign, since   $f(i_k)=f(i_{k+1})$,  and  the sum of these two paths is zero. Hence 
$f_{{\sharp}}\partial(e_{i_0\dots i_n})=0$, and by (\ref{3.6}) diagram (\ref{3.5})
is commutative in the  case  (i).

Consider the case (ii).  Only in this case we are forced to use the asymmetry of $H$. For an allowed path $e_{i_0\dots i_n}\in \mathcal A_n(G)$,  there are no fragments in the form $iji \, (i\ne j)$ in  $e_{f(i_0)\dots f(i_n)}$ since the digraph $H$ is asymmetric. Hence 
\begin{equation}\label{3.8}
\partial f_{\sharp}(e_{i_0\dots i_n})=\partial( e_{f(i_0)\dots f( i_n)}) =\sum_{m=0}^n (-1)^m e_{f(i_0)\dots \widehat{f(i_m)}\dots f(i_n)},
\end{equation}
where  all elementary paths in (\ref{3.8}) are regular. Now we prove  that there are no fragments  $iji \,  (i\ne j)$ in the allowed path  $e_{i_0\dots i_n}\in \mathcal A_n(G)$ in the considered case. If   such a fragment  exists then its  image in the allowed path $e_{f(i_0)\dots f( i_n)}$  has the form  $f(i)f(j)f(i)\,  (f(i)\ne f(j))$   or  
$f(i)f(i)f(i)\,  (f(i)= f(j)) $. The first case is impossible since
digraph $H$ is asymmetric 
and the second case is impossible since in the case (ii) the path $e_{f(i_0)\dots f(i_n)}$ is allowed.  Hence, 
\begin{equation}\label{3.9}
f_{\sharp}\partial(e_{i_0\dots i_n})=f_{\sharp}\left(\sum_{m=0}^n (-1)^m e_{i_0\dots \widehat{i_m}\dots i_n}\right),
\end{equation}
where all elementary paths $e_{i_0\dots \widehat{i_m}\dots i_n}\in \Lambda_n(V_G)$ are regular. Moreover, the image  $f_{\sharp}\left(e_{i_0\dots \widehat{i_m}\dots i_n}\right)\in \Lambda_n(V_H)$ of every such path is also regular since there are no fragments in the form $iji \, (i\ne j)$ in  $e_{f(i_0)\dots f(i_n)}$ as was proved above. Thus we obtain
\begin{equation}\label{3.10}
f_{\sharp}\partial(e_{i_0\dots i_n})=f_{\sharp}\left(\sum_{m=0}^n (-1)^m e_{i_0\dots \widehat{i_m}\dots i_n}\right) =\sum_{m=0}^n (-1)^m e_{f(i_0)\dots \widehat{f(i_m)}\dots f(i_n)}.
\end{equation}
It follows now from (\ref{3.9}) and (\ref{3.10}) that diagram (\ref{3.5}) is commutative in the case (ii). 
\end{proof}

\begin{example}\label{e3.4} \rm Consider the digraph $G=(0\to 1 \to 2)$ and the digraph $H=(0\rightleftarrows 1)$. Define the map $f\colon G\to H$ by setting 
$f(0)=f(2)=0$, $f(1)=1$. Let $e_{012}\in \mathcal A_2(G)$. We compute   homomorphisms in (\ref{3.5})  using (\ref{2.1}) and (\ref{3.4}):
$$
f_{\sharp}\partial (e_{012})= f_{\sharp}(e_{12}-e_{02}+e_{01})= f_{\sharp}(e_{12})- f_{\sharp}(e_{02})+ f_{\sharp}(e_{01})=e_{10}+ e_{01}
$$
since $f_{\sharp}(e_{02})=0$ by (\ref{3.4}). Furthermore, 
$$
\partial f_{\sharp} (e_{012})= \partial(e_{010})= e_{10}- e_{00} +e_{01}.
$$
Hence diagram (\ref{3.5}) is not commutative in general case. 
\end{example} 

\begin{proposition}\label{p3.5} Let  $f\colon G\to H$ be a  map of asymmetric digraphs.  Then  the restriction of the map $f_{\sharp}\colon \mathcal A_n(G)\to \mathcal A_n(H)$ to  $\Pi_n(G)$ satisfies the condition 
\begin{equation}\label{3.11}
f_{\sharp}\left(\Pi_n(G)\right)\subset \Pi_n(H).
\end{equation}
Moreover, this map commutes with the differential $\partial$ and hence, 
defines a morphism of chain complexes $\Pi_*(G) \to  \Pi_*(H)$.
\end{proposition} 
\begin{proof}  Follows from  Lemma \ref{l3.3} and \cite[Lemma 2.3]{Mihyper}. For  $n\geq 0$,   consider  a diagram  
\begin{equation}\label{3.12}
\begin{matrix}
\Pi_{n}(G) & &
 \Pi_{n}(H) &   &\\
 \tau_G \downarrow \ && \ \ \downarrow \tau_H && \ \ \ \ \\
 \mathcal A_{n}(G) &\overset{f_{\sharp}}\longrightarrow&
 \mathcal A_{n}(H)&   &\\
\ \partial \downarrow \ &&\  \ \downarrow\partial && \ \ \ \ \\
\Lambda_{n-1}(V_G) &\overset{f_{\sharp}}\longrightarrow&
 \Lambda_{n-1}(V_H)&& \\
\nu_G\uparrow && \ \ \ \uparrow \nu_H  && \ \ \ \ \\
\mathcal A_{n-1}(G) &\overset{f_{\sharp}}\longrightarrow&
 \mathcal A_{n-1}(H)&   &\\
\end{matrix} 
\end{equation}
where $\nu_G$, $\nu_H$,  $\tau_G$,  and $\tau_H$ are natural inclusions. The diagram  (\ref{3.12}) is commutative by Lemma \ref{3.3} and by the definition of  $f_{\sharp}$.  Let $w\in \Pi_n(G)$.  
To prove  the relation (\ref{3.11})  it is necessary to check that in (\ref{3.12})
$$
\partial \left[ f_{\sharp}\tau_G(w)\right] \in \operatorname{Im} \nu_H.
$$
Since $w\in  \Pi_n(G)$ we have 
\begin{equation}\label{3.13}
\partial \left( \tau_G(w)\right) =\nu_G(v) \ \text{where} \ v\in \mathcal A_{n-1}(G).  
\end{equation}
Furthermore, using relation (\ref{3.13}) and commutativity of diagram (\ref{3.12}) we obtain 
\begin{equation*}
\partial  f_{\sharp}[\tau_G(w)] =f_{\sharp} \partial [\tau_G(w)]= f_{\sharp} 
\nu_G(v)=\nu_Hf_{\sharp}(v)
\end{equation*}
and relation (\ref{3.11}) is proved. Now consideration of commutative upper 
square of (\ref{3.12}) finishes the proof.  
\end{proof}

\begin{corollary}\label{c3.6} The primitive  path homology  is functorial on  the subcategory $\mathcal AD$  of the category $\mathcal D$. 
\end{corollary} 

\begin{theorem}\label{t3.7} Let $G$ be an asymmetric  digraph. 
Then the chain complexes $\Pi_*(G)$ and $\Omega_*(G)$  are naturally 
isomorphic.   
\end{theorem} 
\begin{proof} For $n\geq 0$, the basis of the  module $\mathcal A_n(G)$ is given by the allowed elementary paths $e_{i_0\dots i_n}$ in $G$. 
  Consider the following  commutative diagram 
\begin{equation}\label{3.14}
\begin{matrix}
\mathcal A_{n}(G) &\overset{=}\longrightarrow&
 \mathcal A_{n}(G)&   &\\
\downarrow\ &&  \downarrow\eta_n && \ \ \ \ \\
\Lambda_n(V_G) &\overset{p_n}\longrightarrow& \mathcal R_n(V_G)&= &\Lambda_n(V_G)/I_n(V_G)\\
\partial \downarrow \ &&\downarrow\partial \ && \ \ \ \ \\
\Lambda_{n-1}(V_G) &\overset{p_{n-1}}\longrightarrow&
 \mathcal R_{n-1}(V_G)&=& \Lambda_{n-1}(V_G)/I_{n-1}(V_G)\\
 \ \ \uparrow  && \ \ \ \uparrow \eta_{n-1}&& \ \ \ \ \\
\mathcal A_{n-1}(G) &\overset{=}\longrightarrow&
 \mathcal A_{n-1}(G)&   &\\
\end{matrix} 
\end{equation}
in which  horizontal homomorphisms $p_n, p_{n-1}$   are natural projections,   the right 
differential  is induced by the left  differential,    the upper and bottom left vertical maps are natural inclusions,  and the map $\eta_n$ is given on basic elements $e_{i_0\dots i_n}\in \mathcal A_n(G)$   by  
$\eta_n(e_{i_0\dots i_n}) = [e_{i_0\dots i_n}]$ where 
$[e_{i_0\dots i_n}]\in \mathcal R_n(V_G)$ is the equivalence class defined in Section \ref{S2}.  We have  $p_n(e_{i_0\dots i_n}) = [e_{i_0\dots i_n}]$  for $n\geq 0$ and $p_{-1}=0$. 
The homomorphism $\eta_{n}$  is an inclusion   since $\mathcal A_n(G)\cap I_n(G)=0$. Thus we can naturally identify  $\mathcal A_n(G)$ with the  submodule  of $\mathcal R_n(V_G)$ which is generated by classes $[e_{i_0\dots i_n}]$ for elementary allowed paths $e_{i_0\dots i_n}$. The inclusion 
$\eta_{n-1}$ has the similar properties as $\eta_n$.

For every elementary allowed path $e_{i_0}\in \mathcal A_0(G)$ we have $\partial (e_{i_0})=0\in \mathcal A_{-1}(G)=0$ and $\partial([ e_{i_0}]) =0 \in \mathcal R_{-1}(V_G)=0$. Hence 
\begin{equation}\label{3.15}
 \Pi_0(G)= \Omega_0(G)=\mathcal A_0(G)
\end{equation}
and  $p_0$ gives an isomorphism  $p_0\colon \Pi_0(G)\to \Omega_0(G)$. 

For every elementary allowed path $e_{i_0i_1}\in \mathcal A_1(G)$ we have $\partial (e_{i_0i_1})=e_{i_1}-e_{i_0}\in \mathcal A_{0}(G)$ and 
$\partial([ e_{i_0i_1}]) = [e_{i_1}]-[e_{i_0}]  \in \eta_0(\mathcal A_{0}(G))$. Hence 
\begin{equation}\label{3.16}
 \Pi_1(G)= \Omega_1(G)=\mathcal A_1(G)
\end{equation}
and   $p_1$ gives an isomorphism  
$p_1\colon \Pi_1(G)\to  \Omega_1(G)$. 

For $n\geq 2$,  consider an  elementary  allowed path $e_{i_0\dots i_n}\in \mathcal A_n(G)$.   The sequence of vertices $i_0\dots i_n$ doesn't contain fragments  in the form $iji\,  (i,j \in V_G)$ since there are no symmetric arrows in $G$. Hence,  
each summand  in (\ref{3.17}) below 
\begin{equation}\label{3.17}
\partial(e_{i_{0}...i_{n}})=\sum_{m=0}^n(-1)^m e_{i_{0}\dots \widehat{i_{m}}\dots i_{n}}
\end{equation}
 doesn't belong to $I_{n-1}(V_G)$. Hence 
\begin{equation}\label{3.18}
p_{n-1}  \partial(e_{i_{0}...i_{n}})= \sum_{m=0}^n(-1)^m [e_{i_{0}\dots \widehat{i_{m}}\dots i_{n}}]
\end{equation} 
 Consider an  element 
\begin{equation}\label{3.19}
w=\sum r^{\alpha}e_{\alpha}\in \mathcal A_n(G)\subset \Lambda_n(V_G),  \ 
\  
\partial(w)=\sum r^{\beta}e_{\beta}  \in  \Lambda_{n-1}(V_G) 
\end{equation} 
where $r^{\alpha}, r^{\beta}\in R$,  $e_{\alpha}$ are elementary allowed $n$-paths,  and 
$e_{\beta}$ are elementary $(n-1)$-paths  obtained by the application of $\partial$.  It follows from  diagram (\ref{3.14}),     (\ref{3.17}), and ({3.18})  that the diagram 
\begin{equation}\label{3.20}
\begin{matrix}
\sum r^{\alpha}e_{\alpha} &\overset{p_n} \longrightarrow & \sum r^{\alpha}[e_{\alpha}]\\
\\
\partial \downarrow && \downarrow \partial\\
\\
\sum r^{\beta}e_{\beta} &\overset{p_{n-1}}\longrightarrow & \sum r^{\beta}[e_{\beta}]\\
\end{matrix}
\end{equation} 
is commutative. Similarly to the inclusion $\eta_n$, the inclusion 
$\eta_{n-1}$ is given on basic elements $e_{i_0\dots i_{n-1}}\in \mathcal A_{n-1}(G)$   by  
$\eta_{n-1}(e_{i_0\dots i_{n-1}}) = [e_{i_0\dots i_{n-1}}]$. Hence the commutativity of the bottom square in (\ref{3.14}) and (\ref{3.20}) imply that 
$$
\sum r^{\beta}e_{\beta}\in \mathcal A_{n-1}(G)
$$
 if and only if 
$$
\sum r^{\beta}[e_{\beta}]\in \eta_{n-1}\left(\mathcal A_{n-1}(G)\right)
$$
Hence, for $n\geq -1$, the  restriction  of the homomorphism $p_n$ to the submodule 
$\Pi_n(G)\subset \Lambda_n(V_G)$  gives the natural isomorphism of modules 
$\Pi_{n}(G)\to \Omega_n(G)$ which agrees with the differential as follows from diagram (\ref{3.14}).  
\end{proof}

\begin{corollary}\label{c3.8} Let $G$ be an asymmetric  digraph. 
Then we have a natural isomorphism  $\mathcal H_*(G)\cong H_*(G)$ between  primitive  path  homology groups  and path homology groups. 
\end{corollary}

\begin{proposition}\label{p3.9} Let  $G$ be a digraph and $w \in \mathcal A_n(G)\, (n\geq 0)$.  Then $w\in \Pi_n(G)$ if and only if 
\begin{equation}\label{3.21}
\partial^m(w)\in \mathcal A_{n-1}(G) \ \ \text{for} \ \ 0\leq m\leq n
\end{equation}
where $\partial^m$ is defined in (\ref{2.2}). 
\end{proposition}
\begin{proof} Let $w \in \mathcal A_n(G)$. For $n=0$ it is nothing to prove. Consider the case $n\geq 1$. By  (\ref{2.1}) and (\ref{2.2}),   we have 
\begin{equation}\label{3.22}
\partial(w)=\partial^0(w) + (-1)\partial^1(w)+\dots  +(-1)^n\partial^n(w). 
\end{equation}
If the condition \ref{3.21} is satisfied then every summand in  the   sum  (\ref{3.22}) lays 
in $\mathcal A_{n-1}(G)$ and, hence, $\partial(w)\in \mathcal A_{n-1}(G)$.  By (\ref{3.1}) this means that $w\in \Pi_n(G)$. 

Now we prove the inverse statement.   Let $w\in \Pi_n(G)$ that is $w\in \mathcal A_n(G)$  and   the sum in (\ref{3.22}) lays in $\mathcal A_{n-1}(G)$. It is necessary to prove  that  each element  $\partial^m(w)$ in   (\ref{3.22}) lays in $\mathcal A_{n-1}(G)$.   We can write down   $w$ in the form 
\begin{equation}\label{3.23}
w=\sum_{j=1}^k r^jp_{j}=\sum_{j=1}^kr^je_{i_0^j\dots i_{m-1}^ji_m^ji_{m+1}^j\dots i_n^j}
\end{equation}
where 
 $
p_j=e_{i_0^j\dots i_{m-1}^ji_m^ji_{m+1}^j\dots i_n^j}
$
are pairwise distinct  allowed elementary paths in  $G$ and $r^j\in R$.  By   (\ref{2.1}), for $m=0$ and $m=n$ we have 
\begin{equation}\label{3.24}
\partial^0(w)= \sum_{j=1}^kr^j e_{\widehat{ i_0^j} i_1^j\dots i_n^j}= \sum_{j=1}^kr^j e_{ i_1^j\dots i_n^j}
\in \mathcal A_{n-1}(G)
\end{equation}
and 
\begin{equation}\label{3.25}
\partial^n(w)= \sum_{j=1}^kr^j e_{i_0^j \dots i_{n-1}^j\widehat{ i_n^j}}= \sum_{j=1}^kr^j e_{ i_0^j\dots i_{n-1}^j}
\in \mathcal A_{n-1}(G)
\end{equation}
since each path  $e_{i_0^j\dots  i_n^j}$ 
is an allowed  path in $G$. 

Now consider   $m$ such that  $1\leq m\leq n-1$.    Let $P_m=\{(i,i^{\prime})\}$ be a set of pairs of vertices 
$i, i^{\prime}\in V_G$ 
such that there exists an elementary path  $p_j\, (1\leq j\leq k)$ in (\ref{3.23})  which 
has the form  
 $$
p_j=e_{i_0^j\dots i_{m-1}^ji_m^ji_{m+1}^j\dots i_n^j}  \ \text{where} \ \ 
(i_{m-1}^j, i_{m+1}^j)=(i, i^{\prime}).
$$
Then 
\begin{equation}\label{3.26}
\partial^m(w)=\sum_{(i_{m-1}^j, i_{m+1}^j)\in P_m}r^j e_{i_0^j\dots i_{m-1}^ji_{m+1}^j\dots i_n^j}\in \Lambda_{n-1}(V_G).
\end{equation}
We note the following.  If  $(i_{m-1}^j\to i_{m+1}^j)\in E_G$  then  path $e_{i_0^j\dots i_{m-1}^ji_{m+1}^j\dots i_n^j}$ in (\ref{3.26}) is allowed and lays in $\mathcal A_{n-1}(G)$. If $(i_{m-1}\to i_{m+1})\notin E_G$  then  $e_{i_0^j\dots i_{m-1}^ji_{m+1}^j\dots i_n^j}\notin \mathcal A_{n-1}(G)$.    Let 
$$
T_m=\{(i, i^{\prime})\in P_m \, | \,   (i\to i^{\prime})\in E_G\}
$$
and 
$$
S_m=\{(i, i^{\prime})\in P_m \, |\, (i\to i^{\prime})\notin E_G]\}. 
$$
Thus we can write down $P_m$ in the form 
\begin{equation}\label{3.27}
P_m=T_m\cup S_m \  \text{where } \ T_m\cap S_m=\emptyset.
\end{equation}
Hence,  we can write down  (\ref{3.26})   in the form 
\begin{equation}\label{3.28}
\partial^m(w)=\sum_{(i_{m-1}^j, i_{m+1}^j)\in T_m}r^j\partial^m(p_j)+\sum_{(i_{m-1}^j, i_{m+1}^j)\in S_m}r^j\partial^m(p_j).
\end{equation}
where the first summand lays in $\mathcal A_{n-1}(G)$. 

Since  $\partial(w)\in \mathcal A_{n-1}(G)$ 
 we obtain  using (\ref{3.24}), (\ref{3.25}), and (\ref{3.28}) 
\begin{equation*}
\sum_{1\leq m\leq n-1}\left(     \sum_{(i_{m-1}^j, i_{m+1}^j)\in T_m}r^j\partial^m(p_j)+\sum_{(i_{m-1}^j, i_{m+1}^j)\in S_m}r^j\partial^m(p_j)\right)\in \mathcal A_{n-1}(G).
\end{equation*}
Since all summands by $T_m$ lays in $\mathcal A_{n-1}(G)$ we conclude that 
\begin{equation}\label{3.29}
\sum_{1\leq m\leq n-1}\left(\sum_{(i_{m-1}^j, i_{m+1}^j)\in S_m}r^j\partial^m(p_j)\right)\in \mathcal A_{n-1}(G). 
\end{equation}
Denote by 
\begin{equation}\label{3.30}
\mathcal S^m_{n-1}(G)=\langle \partial^m(p_j)\, | \, (i_{m-1}^j, i_{m+1}^j)\in S_m\rangle \subset \Lambda_{n-1}(V_G)
\end{equation}
 a submodule 
 generated by the elements $\partial^m(p_j)$ such that  $(i_{m-1}^j, i_{m+1}^j)\in S_m$.
We state now, that 
\begin{equation}\label{3.31}
\mathcal S^m_{n-1}(G)\cap \mathcal S^l_{n-1}(G)=0\in \Lambda_{n-1}(V_G) \ \ \text{for} \  \ \ m\ne l.
\end{equation}
For every basic element 
$$
\partial^m(p^j)=e_{i_0^j\dots i_{m-1}^j  i_{m+1}^j\dots i_{n}^j}\in \mathcal S_{n-1}^m(G)\, (1\leq m\leq n-1), 
$$  
there is an arrow $i_t^j\to i_{t+1}^j$  for all $t$ except $t=m-1$. Similarly,  for every basic element 
$$
\partial^l(p^j)=e_{i_0^j\dots   i_{l-1}^j  i_{l+1}^j\dots i_{n}^j}\in \mathcal S_{n-1}^l(G)\, (1\leq l\leq n-1),
$$ 
there is an arrow $i_t^j\to i_{t+1}^j$  for all $t$ except $t=l-1$. Hence,  the condition (\ref{3.31}) is satisfied. Recall that by the definition of the set $\mathcal S^m_{n-1}$,   for every element $\partial^m(p^j)=e_{i_0^j\dots i_{n}^j}\in \mathcal S_{n-1}(G) \subset \Lambda_{n-1}(V_G)$ we have 
\begin{equation}\label{3.32}
\partial^m(p^j)=e_{i_0^j\dots i_{n}^j}\notin \mathcal A_{n-1}(G) . 
\end{equation}
Hence, by (\ref{3.29}), (\ref{3.31}), and (\ref{3.32}) we obtain 
 that for $1\leq m\leq n-1$ 
\begin{equation}\label{3.33}
\sum_{(i_{m-1}^j, i_{m+1}^j)\in S_m}\partial^m(p_j)=0\in \mathcal A_{n-1}(G). 
\end{equation}
By the consideration above,  the first summand in (\ref{3.28}) lays in $\mathcal A_{n-1}(G)$  and  by (\ref{3.33})  the second  summand in (\ref{3.28}) is zero.   Hence $\partial^m_n(w)\in \mathcal A_{n-1}(G)$ and the proposition is proved. 
\end{proof}

\begin{definition}\label{d3.10}  \rm  Let  $G=(V,E)$ be a digraph.   We denote by  $\overline{G}=(\overline{V}, \overline E)$ a digraph with 
$$
\overline V=V \ \  \text{and} \ \  \overline{E}=\{v\to w\, | \, v, w\in V,  (w\to v)\in E\}.
$$ 
We call the digraph $\overline{G}$ by \emph{inverse directed digraph} of the digraph $G$.   
\end{definition} 

Thus,  the digraph  $\overline{G}$  is obtained from $G$ by changing directions of all arrows.
It easy to see that  the operation $G\to \overline{G}$ is an involution on the set of digraphs. 

\begin{theorem}\label{t3.11} For every digraph $G=(V,E)$,  we have an isomorphism of chain complexes 
$$
\Pi_*(G) \to \Pi_*(\overline{G})
$$ 
which induces an isomorphism of primitive path homology groups 
$$
\mathcal  H_n(G) \to 
\mathcal H_n( \overline{G})  \ \ \text{for} \ \ n\geq 0.
$$ 
\end{theorem}
\begin{proof}  For the convenience of  readers  we give the  proof which  is similar to \cite[Theorem 4.14]{Mi2}. For $n\geq 0$, define an isomorphism of  $R$-modules 
$$
\tau\colon \Lambda_n(V) \to \Lambda_n(V)  \  \ \text{by}   \ \ 
\tau(e_{i_0\dots i_n})=\begin{cases} e_{i_n\dots i_0} & \text{for}\  n=0,3\bmod 4, \\
e_{i_n\dots i_0} & \text{for}\  n=1,2\bmod 4. 
\end{cases} 
$$
For $n=-1$, we set $\tau=0\colon 0=\Lambda_{-1}(V)\to \Lambda_{-1}(V)=0$. 
The isomorphism $\tau$ commutes with the differential $\partial$ in (\ref{2.1}) and, hence, 
defines  an isomorphism of chain complexes $\tau_*\colon \Lambda_*(V)\to \Lambda_*(V)$. 
For $n\geq -1$, let  $\mathcal A_n(G)\subset \Lambda_n(V)$ be a submodule generated by allowed elementary paths.  By Definition \ref{d3.10} the restriction of the isomorphism $\tau$ 
defines an isomorphism $\tau\colon \mathcal A_n(G)\to  \mathcal A_n(\overline G)$.  Let
$w\in \Pi_n(G)$,  that is $w \in \mathcal A_{n}(G)$    and $\partial w \in \mathcal A_{n-1}(G)$.  
Then $\tau(w)\in  \mathcal A_n(\overline G)$,  $\partial \tau(w) =\tau \partial(w)\in   \mathcal A_{n-1}(\overline G)$ and, hence,   $\tau(w)\in \Pi_n(\overline G)$. The inverse conclusion  is similar since $\tau_*$ and the restriction of $\tau$  are isomorphisms. 
\end{proof}

\section{ Cluster digraphs}\label{S4}
\setcounter{equation}{0}

In this Section we  define a cluster  subgraph,  a tail subgraph, and a head subgraph of a digraph $G$ and study primitive  path homology of these digraphs.   We describe relations between these homology  and their relation to  primitive path homology of the ambient digraph. We fix an unitary and commutative  ring $R$ of coefficients.  All modules are $R$-modules and we don't mention $R$ in  notations. 

\begin{definition}\label{d4.1} \rm   For $a,b\in V_G$, define $(a,b)$-\emph{cluster digraph} $G^{[a,b]}=(V^{[a,b]}, E^{[a,b]})$  of a digraph $G$  as a minimal  subgraph  which contains   all vertices and all  arrows  of every  $(a,b)$-cluster path $w^{[a,b]}\in \mathcal A_n(G)$ for $n=0,1,\dots $.
\end{definition} 

Suppose a cluster digraph $G^{[a,b]}$ exists. Then $\mathcal A_n(G^{[a,b]})\subset \mathcal A_{n}(G)$ for $n\geq -1$ by  Definition \ref{d4.1}.  For every vertex $a\in V_G$ the cluster digraph  $G^{[a,a]}$ exists and contains the vertex $a$.

We note also that a cluster digraph $G^{[a,b]}$ may be not defined  for some pairs $(a,b)$  of vertices. For example, for the digraph $G=(0\to 1)$ the cluster digraph $G^{[1,0]}$ is  not  defined.

\begin{proposition}\label{p4.2}  Let  $G^{[a,b]}$ be a  cluster digraph in $G$. Then it is an induced subgraph.   
\end{proposition} 
\begin{proof}   Let $k$ be a vertex of an allowed path $e_{a\dots sk\dots b}\in \mathcal A_n\left(G^{[a,b]}\right)$ and $l\ne k $ be a vertex of an allowed path $e_{a\dots lm\dots b}\in \mathcal A_n\left(G^{[a,b]}\right)$  such that  $(k\to l)\in E_G$. Then the allowed  path $e_{a\dots sklm\dots b}$ from $a$ to $b$
 contains the arrow $k\to l$ and, hence, 
$(k\to l)\in E^{[a,b]}$. 
\end{proof}

Let  $G$ be  a digraph and  $w\in\mathcal A_n(G)$ be an  allowed path.  For $n\geq 0$,  we can write 
$w$ in the form  
\begin{equation}\label{4.1} 
w=\sum_{a\in V}w^{[a, \cdot]}   \ \text{where}  \
w^{[a, \cdot]} =\sum_{ai_1\dots i_n}r^{ai_1\dots i_n}e_{ai_1\dots i_n}, \ \ r^{ai_1\dots i_n}\in R
\end{equation}
and  $e_{ai_1\dots i_n}$ are elementary allowed paths and 
$
a=t\left(e_{ai_1\dots i_{n-1}i_n}\right). 
$
We call the   path   $w^{[a, \cdot]}$ by \emph{$a$-tail  path}. 
Similarly, for $n\geq 0$,  we can write 
$w$ in the form  
\begin{equation}\label{4.2} 
w=\sum_{b\in V}w^{[\cdot,b]},    \ \  \  
w^{[\cdot,b]} =\sum_{i_0\dots i_{n-1}b}r^{i_0\dots i_{n-1}b}e_{i_0\dots i_{n-1}b}, \ \ r^{i_0\dots i_{n-1}b}\in R
\end{equation}
and $e_{i_0\dots i_{n-1}b}$ are elementary allowed paths and
$
b=h\left(e_{i_0i_1\dots i_{n-1}b}\right). 
$
We call the   path  $w^{[\cdot,b]}$ by \emph{$b$-head  path}.

\begin{theorem}\label{t4.3} Let  $G$ be a digraph and  $ \Pi_n(G)\, (n\geq 0)$ be a module in (\ref{3.1}). 

{\rm (i)} A path  $w \in \mathcal A_n(G)$ lays in $\Pi_n(G)$ if and only  if $w^{[a, \cdot]}\in \Pi_n(G)$ for every tail fixed path $w^{[a, \cdot]}$  in (\ref{4.1}).

 {\rm (ii)} A path  $w \in \mathcal A_n(G)$ lays in $\Pi_n(G)$  if and only  if $w^{[\cdot, b]}\in \Pi_n(G)$ for every head fixed path $w^{[\cdot, b]}$  in (\ref{4.2}).

 {\rm (iii)} A path  $w \in \mathcal A_n(G)$ lays in $\Pi_n(G)$  if and only  if $w^{[a, b]}\in \Pi_n(G)$ for every cluster path $w^{[a, b]}$  in (\ref{2.7}).
\end{theorem}
\begin{proof}  For $n=0$ it is nothing to prove. For $n\geq 1$, consider  the case (i).    If  $w^{[a.\cdot]}\in \Pi_n(G)$ for all tail paths in (\ref{4.1}) then
$\partial(w^{[a.\cdot]})\in \mathcal A_{n-1}(G)$ and 
$$
\partial(w)=  \partial\left(\sum_{a\in V}w^{[a.\cdot]}\right)=\sum_{a\in V}\partial\left(w^{[a.\cdot]}\right)\in \mathcal A_{n-1}. 
$$ 
Hence $w\in \Pi_n(G)$. 

Now we prove the inverse statement.    Let $w\in \Pi_n(G)$ and, hence, $\partial w\in \mathcal A_{n-1}(G)$.  The paths $w^{[a.\cdot]}$ in (\ref{4.1}) are allowed and it is necessary to prove that 
$\partial\left(w^{[a.\cdot]}\right)\in \mathcal A_{n-1}(G)$ for every $w^{[a.\cdot]}$. Every  path $w^{[a,\cdot]}$ in the decomposition (\ref{4.1}) is a  linear combination of allowed paths with the tail vertex $a$: 
$$
w^{[a.\cdot]}=\sum_{k} r^{k}e_{ai_1^k\dots i_n^k},  \ \  r^{k}\in R. 
$$
 We have 
\begin{equation}\label{4.3}
\partial\left(w^{[a,\cdot]}\right)=\sum_{k} r^{k}e_{i_1^ki_2^k\dots i_n^k} + 
\sum_k\left(\sum_{j=1}^n  r^{k}e_{ai_1^k\dots \, \widehat{i_j^k}\,\dots i_n^k}\right). 
\end{equation}
For every $a\in V$, the first sum in (\ref{4.3}) is an allowed path.  The second sum can contains allowed and non-allowed elementary paths. An non-allowed elementary path in the second sum  can not be cancel  with a similar  non-allowed elementary fitting into $w^{[c, \cdot]}$ with $c\ne a$. 
Since there are no non-allowed paths in $\partial w$ we conclude that they have to cancel in the second sum in (\ref{3.4}). Hence $\partial \left(w^{[a.\cdot]}\right)\in \mathcal A_{n-1}(G)$ and  $w^{[a.\cdot]}\in \Pi_{n}(G)$. The proof in the  cases (ii) and (iii) is  similar. 
\end{proof} 

\begin{corollary}\label{c4.4} {\rm (i)} Let  $G=(V,E)$ be a digraph and $n\geq 0$.  The module $\Pi_n(G)$  is a direct sum of submodules 
$$
\Pi_n(G)=\bigoplus_{a\in V }  \Pi_n^{[a,\cdot]}(G)$$
where $ 
\Pi_n^{[a,\cdot]}(G)=\left\langle w^{[a,\cdot]}\in \Pi_n(G) \right\rangle
$
is  generated by 
$a$-tail  paths  in $\Pi_n(G)$.  If there are no $a$-tail paths    in dimension $n$ we set $\Pi_n^{[a,\cdot]}(G)=0$. 

{\rm (ii)}  The module $\Pi_n(G)$  is a direct sum of submodules 
$$
\Pi_n(G)=\bigoplus_{b\in V }  \Pi_n^{[\cdot, b]}(G) $$
where $
\Pi_n^{[\cdot, b]}(G)=\left\langle w^{[\cdot, b]}\in \Pi_n(G) \right\rangle
$
is  generated by 
$b$-head  paths  in $\Pi_n(G)$.  If there are no $b$-head paths   in dimension $n$ we set $\Pi_n^{[\cdot, b]}(G)=0$. 

{\rm (iii)}  The module $\Pi_n(G)$  is a direct sum of submodules 
$$
\Pi_n(G)=\bigoplus_{a, b\in V }  \Pi_n^{[a, b]}(G) $$
where $\Pi_n^{[a, b]}(G)=\left\langle w^{[a, b]}\in \Pi_n(G) \right\rangle
$
is  generated by 
$(a,b)$-cluster  paths  in $\Pi_n(G)$.  If there are no $(a,b)$-cluster paths   in dimension $n$ we set $\Pi_n^{[a, b]}(G)=0$.
\end{corollary}

\begin{definition}\label{d4.5} \rm  (i)  For $a\in V$, we define an 
$a$-\emph{tail digraph} $G^{[a,\cdot]}=(V^{[a,\cdot]}, E^{[a,\cdot]})$  of a digraph $G$  as a minimal  subgraph  which contains   all vertices and all  arrows
  of every  $a$-tail path $w^{[a,\cdot]}\in \mathcal A_n(G)$ for $n=0,1,\dots $. 

(ii)  For $b\in V$, we define an 
$b$-\emph{head digraph} $G^{[\cdot, b]}=(V^{[\cdot, b]}, E^{[\cdot, b]})$  of a digraph $G$  as a minimal  subgraph  which contains   all vertices and all  arrows
  of every  $b$-head path $w^{[\cdot, b]}\in \mathcal A_n(G)$ for $n=0,1,\dots $. 
\end{definition} 

\begin{proposition}\label{p4.6}  Every $a$-tail and every  $b$-head digraph of a digraph $G$  is an induced subgraph.   
\end{proposition} 
\begin{proof}  Similar to the proof of Proposition \ref{p4.2}. 
\end{proof} 

\begin{proposition}\label{p4.7} Let $G=(V,E)$ be a digraph and $n\geq 0$. For $a,b\in V$ we have the following equalities of modules 
\begin{equation}\label{4.4}
\begin{matrix}
 \Pi_n^{[a,\cdot]}(G)=\Pi_n^{[a,\cdot]}\left(G^{[a,\cdot]}\right),  \ \   \Pi_n^{[\cdot, b]}(G)=\Pi_n^{[\cdot, b]}\left(G^{[\cdot, b]}\right),  \\
\\
  \Pi_n^{[a,b]}(G)=\Pi_n^{[a,b]}\left(G^{[a,b]}\right), \ \text{and} \ \
\Pi_n^{[a,b]}(G)=\Pi_n^{[a, \cdot ]}(G) \cap \Pi_n^{[\cdot, b]}(G).
\end{matrix}
\end{equation}
\end{proposition}
\begin{proof} The equalities in (\ref{4.4}) follow from Definitions \ref{d4.1} and \ref{d4.5},  Theorem \ref{t4.3} and Corollary  \ref{c4.4}. 
\end{proof} 

\begin{remark}\label{r4.8} \rm We note that the image of the restriction  
$$
\partial\colon \Pi_n^{[a,b]}(G) \to \Pi_{n-1}(G)
$$
of the differential $\partial\colon \Pi_n(G)\to \Pi_{n-1}(G)$ 
doesn't belong to the  submodule 
$\Pi_{n-1}^{[a,b]}(G)\subset \Pi_n(G)$ and, hence, 
the modules  $\Pi_n^{[a,b]}(G)$ with the differential $\partial$ doesn't  form a chain complex. The similar statement is true for  the  modules
$ \Pi_n^{[a,\cdot]}(G)$  and $\Pi_n^{[\cdot, b]}(G)$ in (\ref{4.4}). 
\end{remark}

\section{Primitive homology of cluster digraphs}\label{S5}
\setcounter{equation}{0}

In this Section, for a digraph $G$, we construct a primitive path  homology theory  for  a  set of $(a,b)$-cluster paths, that is the set of paths which have the  fixed tail  vertex $a$ and  the fixed head vertex $b$.  
The construction of the chain complex  is similar to the construction of the  chain  complex in  the  primitive path homology theory in Section \ref{S3}.  We note that the homology theory of cyclic paths  $e_{i_{0}\dots i_{n-1}i_n}$ with $i_0=i_n$ was constructed and studied  in \cite{Circuits} and  \cite{MurJC}.
We describe relations between introduced  homology  and  the primitive  path homology and  present several examples. 

Fix a unitary and commutative ring $R$ of coefficients.  The next  definition is a natural  generalization  of   Definition \ref{d2.2} in which  a cluster path is defined  for  a digraph $G$.

\begin{definition}\label{d5.1}\rm  Let  $V$ be a finite set of vertices,   $a,b\in V$ be two vertices, and  $n\geq 0$.  A path  in $\Lambda_n(G)$  is  an $(a,b)$-\emph{cluster path} if it is a linear combination of paths with the tail vertex $a$ and with the head vertex $b$.  The cluster path is denoted $w^{[a,b]}$. The elementary paths fitting into $w^{[a,b]}$ are called \emph{elementary}  $(a,b)$-\emph{cluster paths}.  
\end{definition} 

 Let  $n\geq 1$.  Denote by   $\Lambda _{n}^{[a,b]}=\Lambda _{n}^{[a,b]}( V)\subset  \Lambda_n(V)$  a free submodule generated by    elementary   $(a,b)$-\emph{cluster path}.  The module  $\Lambda _{n}^{[a,b]}$ is called $(a,b)$-cluster module and the elements of $\Lambda_{n}^{[a,b]} $ are called $(a,b)$-\emph{cluster paths   in dimension}  $n$.   We set    also 
$\Lambda_{0}^{[a,b]}=0$.

For $n\geq 1$, define a  \emph{differential} 
\begin{equation}\label{5.1}
d\colon \Lambda _{n}^{[a,b]}\rightarrow \Lambda _{n-1}^{[a,b]}
\end{equation} 
 by setting on     basic elements 
\begin{equation}\label{5.2}
 d (e_{ai_{0}\dots  i_{n-2}b} )=\begin{cases} \sum\limits_{m=0}^{n-2}\left( -1\right)^{m}e_{ai_{0}\dots \widehat{i_{m}}\dots i_{n-2}b}& \text{for} \ n\geq 2,\\
 0& \text{for} \  n=1.
\end{cases}
\end{equation}
 For $n\geq 2$, we can write down the differential $d$  in the form
$$
d (e_{ai_{0}\dots i_{n-2}b})=\sum\limits_{m=0}^{n-2}(-1)^{m} d^m(e_{ai_{0}\dots i_{n-2}b}). 
$$
where $d^m\colon \Lambda _{n}^{[a,b]}\rightarrow \Lambda _{n-1}^{[a,b]}$ is a homomorphism given on the  basic elements by 
\begin{equation*}
d^m(e_{ai_{0}\dots i_{n-2}b})= e_{ai_{0}\dots \widehat{i_{m}}\dots i_{n-2}b}. 
\end{equation*}

Consider the augmented chain complex $\widetilde{\Lambda}_*(V)$  defined in Section \ref{S2}. We have $\widetilde{\Lambda}_n(V)=\Lambda_n(V)$  for 
$n\geq 0$,  $\widetilde{\Lambda}_{-1}(V)=R$,  and 
$\widetilde{\Lambda}_{-2}(V)=0$.

For $n\ge 0$, define an  isomorphism  of modules 
\begin{equation}\label{5.3}
\pi^{[a,b]}\colon \Lambda_n^{[a,b]}(V)\to \widetilde \Lambda_{n-2}(V),  \ \ 
\end{equation}
by setting:
$\pi^{[a,b]}(e_{ai_{0}\dots i_{n-2}b})=e_{i_{0}\dots i_{n-2}}$ for $n\geq 2$,  
$\pi^{[a,b]}(e_{ab})= 1\in R$
 for $n= 1$, and $\pi^{[a,b]}=0$  for $n=0$. 

\begin{proposition}\label{p5.2}  For $n\geq 0$, the following diagram 
\begin{equation}\label{5.4} 
\begin{matrix}
\Lambda_n^{[a,b]}(V)& \overset{\pi^{[a,b]}}\longrightarrow & \widetilde \Lambda_{n-2}(V)\\
\\
d\downarrow&&\downarrow \partial\\
\\
\Lambda_{n-1}^{[a,b]}(V)& \overset{\pi^{[a,b]}}\longrightarrow 
& \widetilde \Lambda_{n-3}(V)\\
\end{matrix}
\end{equation} 
is commutative. 
\end{proposition}
\begin{proof} Follows from the definition of the isomorphism $\pi^{[a,b]}$ and definitions of  differential $\partial$ in (\ref{2.1}) and differential $d$ in (\ref{5.2}). 
\end{proof} 

\begin{corollary}\label{c5.3} For $n\geq 0$, the modules $\Lambda_n^{[a,b]}$ with the differential $d$ form a chain complex.
\end{corollary} 
\begin{proof} In diagram (\ref{5.4})   the modules  
$\widetilde \Lambda_{n}(V)\, (n\geq -2)$ with the differential $\partial$ form a chain complex.  The horizontal maps are isomorphisms. Now the statement follows from commutativity of (\ref{5.4}). 
\end{proof}

\begin{definition}\label{d5.4} \rm Consider  a digraph $G=(V, E)$ and   two vertices $a, b\in V$. For $n\geq 0$, 
let  $\mathcal A_n^{[a,b]}=\mathcal A_n^{[a,b]}(G)$ be  a module generated 
by elementary $(a,b)$-cluster paths  in $G$. If there are no such paths  we set 
$\mathcal A_n^{[a,b]}=0$.  We call the module  $\mathcal A_n^{[a,b]}$ by an  \emph{allowed  $(a,b)$-cluster module}. 
\end{definition} 
It follows from Definition \ref{d5.4} that  for every digraph $G$ we have 
\begin{equation}\label{5.5}
\begin{matrix}
\mathcal A_0^{[a,b]}(G)=0 \ \text{for all} \  a, b\in V,\\   
  \mathcal A_1^{[a,b]}(G)=\begin{cases} \langle e_{ab} \rangle, &  \text{if} \  a\to b \in E, \\
                                                                0,&  \text{otherwise}. \\
\end{cases}\\
 \end{matrix}
\end{equation}

For $n\geq 0$, the module $\mathcal A_n^{[a,b]}(G)$ is a submodule of $\Lambda_{n}(V)$. For $n\geq 1$, define a submodule
\begin{equation}\label{5.6}
\Theta_n^{[a,b]}=\Theta_n^{[a,b]}(G)=\{w\in \mathcal A_n^{[a,b]}(G)\, |
\, d(w)\in  \mathcal A_{n-1}^{[a,b]}(G)\},
\end{equation} 
where $d\colon  \Lambda_{n}^{[a,b]}(V)\to  \Lambda_{n-1}^{[a,b]}(V)$ is the differential in 
(\ref{5.2}). We set $\Theta_0^{[a,b]}=0$. 

\begin{proposition}\label{p5.5} The modules $\Theta_n^{[a,b]}(G)$ with the restriction of the differential $d$ form  a  chain complex $\Theta_*^{[a,b]}=\Theta_*^{[a,b]}(G)$.
\end{proposition} 
\begin{proof}
Similarly to the case of modules $\Pi_n(G)$ defined in (\ref{3.1}), we obtain
that $d(\Theta_n^{[a,b]})\subset \Theta_{n-1}^{[a,b]}$ for $n\geq 1$ and the result follows.
\end{proof} 

 By  (\ref{5.5}) and (\ref{5.6},   for any digraph $G=(V,E)$ we have 
\begin{equation}\label{5.7}
\Theta_1^{[a,b]}(G)=\begin{cases} \langle e_{ab} \ \rangle, &  \text{if} \  a\to b \in E, \\
                                                                0,&  \text{otherwise}. \\
\end{cases}
\end{equation} 

\begin{definition}\label{d5.6}  \rm The homology groups of the chain complex 
$\Theta_*^{[a,b]}(G)$ are called \emph{primitive $(a,b)$-cluster homology groups} of the digraph $G$ and are denoted by $\mathcal H_n^{[a,b]}(G)=\mathcal H_n^{[a,b]}(G,R)\, (n\geq 0)$.
\end{definition} 

Let  $I=(0\to 1)$ be the digraph.   For $n\geq 0$, we define \emph{$n$-cube digraph} $I^n$ by setting $I^{0}=\{0\}$ and
\begin{equation*}
I^{n}=\underset{n\ \text{times}}{\underbrace{I\Box I\Box I\Box \dots \Box I}}
\ \  \text{for}\ \  n\geq 1.
\end{equation*}

\begin{example}\label{e5.7} \rm Let $\mathbb Z$ be the ring of coefficients. 

(i) Consider a $2$-cube digraph $I^2$ below 
$$
\begin{matrix}
1 & \to &3\\
\uparrow && \uparrow \\
0&\to &2.\\
\end{matrix}
$$
We have  
$
\mathcal A^{[0,3]}_0=\mathcal A^{[0,3]}_1=0$, $\mathcal A^{[0,3]}_2=\langle e_{013}, e_{023}\rangle \cong \mathbb Z^2$, $\mathcal A^{[0,3]}_n=0$  for  $n\geq 2$. By the definition of  $d$ in (\ref{5.2}),  we obtain that $d(e_{013}-e_{023})= e_{03}-e_{03}=0$.  Hence, 
$\Theta_2^{[0,3]}=\langle e_{013}-e_{023}\rangle \cong \mathbb Z$ and $\Theta_n^{[0,3]}=0$ for $n\ne 2$. Thus we conclude that 
$
\mathcal H_2^{[0,3]}(I^2)= \mathbb Z$ for  $n=2$ and $\mathcal H_n^{[0,3]}(I^2)=0$ otherwise. 

(ii)  Let $I^3$  be a  digraph cube in Figure \ref{Q}. 
 \begin{figure}[H]
\centering
\begin{tikzpicture}[scale=0.85]
\node (1) at (4,3) {$0$};
\node (2) at (6,3) {$1$};
\draw (1) edge[ color=black!120, thick,-> ] (2);
\node (1a) at (4,5) {$4$};
\node (2a) at (6,5) {$5$};
\node (3a) at (7.4,4.1) {$3$};
\node (30a) at (5.4,4.1) {$2$};
\draw (3a) edge[ color=black!120, thick,<- ] (30a);
\node (3b) at (7.4,6.1) {$7$};
\node (30b) at (5.4,6.1) {$6$};
\draw (3b) edge[ color=black!120, thick,<- ] (30b);
\draw (3a) edge[ color=black!120, thick, ->] (3b);
\draw (30a) edge[ color=black!120, thick, ->] (30b);
\draw (1) edge[ color=black!120, thick, ->] (30a);
\draw (1a) edge[ color=black!120, thick,-> ] (2a);
\draw (2a) edge[ color=black!120, thick,<-] (2);
\draw (1a) edge[ color=black!120, thick,<-] (1);
\draw (3a) edge[ color=black!120, thick, <-] (2);
\draw (1a) edge[ color=black!120, thick,-> ] (30b);
\draw (2a) edge[ color=black!120, thick,-> ] (3b);
\end{tikzpicture}  
\caption{Digraph cube $I^3$.}
\label{Q}
\end{figure}
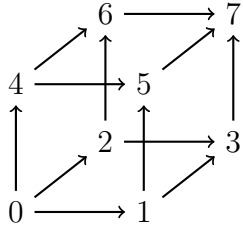
 We have 
$
\mathcal A^{[0,7]}_i=0
$ for $i\ne 3$ and 
$
\mathcal A^{[0,3]}_3=
\langle e_{0137},  e_{0237},  e_{0157}, e_{0457},
 e_{0467}, e_{0267}\rangle.  
$
We have the following differential $d\colon \mathcal  A^{[0,7]}_3\to \Lambda_2^{[0,7]}$\ :
\begin{equation}\label{5.8}
\begin{matrix}
d(e_{0137})&=& e_{037}-e_{017}, &
d(e_{0237})&=& e_{037}-e_{027}, \\
  d(e_{0157})&=&  e_{057}-e_{017}, &
d(e_{0457})&= &e_{057}-e_{047},\\
d( e_{0467})&= &e_{067}-e_{047},&
d(e_{0267})&=&e_{067}-e_{027}. \\
\end{matrix}
\end{equation} 
We can check directly that the kernel of  $d$ is generated by the element 
$$
w=e_{0137}-e_{0237}-e_{0157}+e_{0457}-e_{0467}+e_{0267}. 
$$
Hence, $\Theta_3^{[0,7]}=\langle w\rangle \cong \mathbb Z$  and $\Theta_n^{[0,7]}=0$ for $n\ne 3$. Thus
$
\mathcal H_3^{[0,7]}(I^3)= \mathbb Z$ and $\mathcal H_n^{[0,7]}(I^3)=0$ otherwise.

(iii)  Now we consider a digraph $G$ in Figure \ref{G} below. 

\begin{figure}[H]
\bigskip\bigskip \bigskip \setlength{\unitlength}{0.14in} \centering
\begin{picture}(16,16)(2,2)
\centering
\thicklines
 \put(5,5){\vector(1,0){16}}
\put(4.2,4){$0$}
 \put(5,5){\vector(1,1){10.1}}

\put(21.2,4){$1$}
 \put(21,5){\vector(0,1){16}}
\put(21.3,21.1){$7$}
 \put(5,5){\vector(0,1){16}}
\put(4, 20.9){$4$}
\put(5,21){\vector(1,0){16}}
\put(14.9, 16){$8$}
\put(10.3, 16){$6$}
\put(5,21){\vector(1,-1){5.4}}
 \put(10.5,15.5){\vector(1,0){5}}
 \put(15.5,10.5){\vector(0,1){5}}
\put(21,21){\vector(-1,-1){5.4}}
\put(21,5){\vector(-1,1){5.4}}
\put(16.1, 10.5){$5$}
 \put(13.2,10.5){\vector(1,0){2.4}}
\put(10.5,13.2){\vector(0,1){2.4}}
\put(5,5){\vector(2,3){5.5}}
\put(5,5){\vector(3,2){8.2}}
\put(12.4,9.1){$2$}
\put(9,12.6){$3$}
 \put(13,10.5){\vector(1,2){2.4}}
\put(10.6,13.1){\vector(2,1){4.5}}
\end{picture}
\caption{ The planar digraph $G$ from Example \ref{e5.7} (iii). }
\label{G}
\end{figure}
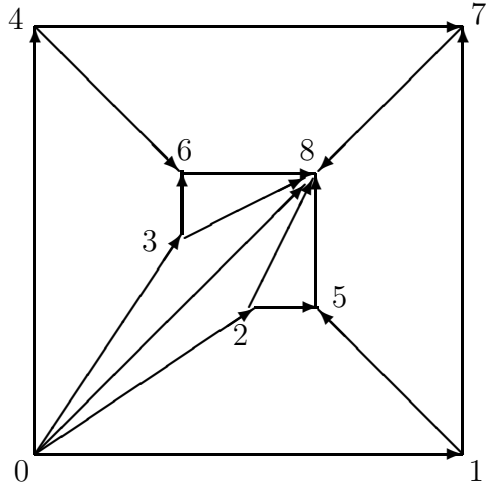

We have 
\begin{equation}\label{5.9}
\mathcal A^{[0,8]}_n=\begin{cases} 0 & \text{for} \  n=0,  \\
\langle e_{08}\rangle\cong \mathbb Z&   \text{for} \  n=1, \\
\langle e_{028}, e_{038}\rangle \cong \mathbb Z^2&  \text{for} \  n=2,  \\
\langle e_{0178},  e_{0158},  e_{0258}, e_{0368},
 e_{0468}, e_{0478}\rangle \cong \mathbb Z^6&  \text{for} \  n=3,  \\
0 & \text{for}  \ n\geq 4.\\
\end{cases}
\end{equation}
We have the following differential  in dimensions 1  and 2: 
\begin{equation}\label{5.10}
d( e_{08})=0, \  d(e_{028})=d(e_{038})=e_{08}. 
\end{equation}
In  dimension 3,  we have the following differential 
$d\colon \mathcal  A^{[0,8]}_3\to \Lambda_2^{[0,7]}$\ :
\begin{equation}\label{5.11}
\begin{matrix}
d(e_{0178})&=& e_{078}-e_{018}, &
d(e_{0158})&=& e_{058}-e_{018}, \\
  d(e_{0258})&=&  e_{058}-e_{028}, &
d(e_{0368})&= &e_{068}-e_{038}, \\
d( e_{0468})&= &e_{068}-e_{048}, &
d(e_{0478})&=&e_{078}-e_{048}. \\
\end{matrix}
\end{equation} 
Any nontrivial linear combination of elements of the right part of equalities  in (\ref{5.11})  doesn't  belong to $\mathcal A_2^{[0, 8]}$ and, hence, by (\ref{5.9}) and (\ref{5.10}) we obtain 
$$
\mathcal H_n^{[0,8]}(G)= \begin{cases} \mathbb Z&  \text{for} \ n=2,\\
                                                                              0&  \text{otherwise}.\\
\end{cases} 
$$
It is interesting to note, the path homology group of the digraph $G$ is trivial in dimension two and the  singular  cubical homology group is non-trivial in dimension two \cite[Propositions 13,14]{Mi4}. 
\end{example}

\begin{theorem}\label{t5.8}  Let $G=(V,E)$ be a digraph and  $a,b\in V$ be two vertices. 

(i)  If $(a\to b)\notin E$ then there is a morphism of chain complexes 
\begin{equation}\label{5.12} 
\begin{matrix}
\Theta_n^{[a,b]}(G)& \overset{\pi^{[a,b]}}\longrightarrow & \Pi_{n-2}(G) \ \ \ (n\geq 1),\\
\end{matrix}
\end{equation} 
where $\Pi_{n}(G)$ is the chain complex  in (\ref{3.1}) for primitive path homology groups and $\pi^{[a,b]}$ is the restriction of the homomorphism $\pi^{[a,b]}$ defined in (\ref{5.3}). 

(ii) If $(a\to b)\in E$ then there is a morphism of chain complexes 
\begin{equation}\label{5.13} 
\begin{matrix}
\Theta_n^{[a,b]}(G)& \overset{\pi^{[a,b]}}\longrightarrow & \widetilde{\Pi}_{n-2}(G)  \ \ \ (n\geq 1),\\
\end{matrix}
\end{equation} 
where $\widetilde{\Pi}_{*}(G)$ is the reduced chain complex  and $\pi^{[a,b]}$ is the restriction of the homomorphism $\pi^{[a,b]}$ defined in (\ref{5.3}). 
\end{theorem}

\begin{proof} (i)  In this case,  we have   $\mathcal A_1^{[a,b]}(G)= \mathcal A_{-1}(G)=0$ and the restriction of the  homomorphism $\pi^{[a,b]}$ to
\begin{equation}\label{5.14}
\pi^{[a,b]}\colon \mathcal A_2^{[a,b]}(G)\to \mathcal A_{0}(G)
\end{equation}
is given by $\pi^{[a,b]}\left(e_{ai_0b}\right)=e_{i_0}\in V$. Hence 
$\Theta_2^{[a,b]}=\mathcal A_2^{[a,b]}$, $\Pi_0=\mathcal A_0$, and 
the homomorphism in (\ref{5.16}) is the   homomorphism 
$\Theta_2^{[a,b]}(G) \to \Pi_0(G)$. Now consider the case 
$n\geq 3$.  For  every elementary  $(a,b)$-cluster  path 
$e_{ai_0\dots i_{n-2}b}\in \mathcal A_n^{[a,b]}(G)$, we have 
$
\pi^{[a,b]}(e_{ai_0\dots i_{n-2}b})= e_{i_0\dots i_{n-2}}\in \mathcal A_{n-2}(G).
$
Consider  a  diagram 
\begin{equation}\label{5.15} 
\begin{matrix}
& &\mathcal A_n^{[a,b]}(G)& \overset{\pi^{[a,b]}}\longrightarrow & \mathcal A_{n-2}(G) & & \\
\\
& &d\downarrow&&\downarrow \partial & & \\
\\
\mathcal A_{n-1}^{[a,b]}(G)&\overset{\subset}\to &\Lambda_{n-1}^{[a,b]}(V)& \overset{\pi^{[a,b]}}\longrightarrow 
&  \Lambda_{n-3}(V)&\overset{\supset}\leftarrow & \mathcal A_{n-3}(G)\\
\end{matrix}
\end{equation} 
in which   left and right bottom horizontal  maps are the natural inclusions. Diagram (\ref{5.17}) is commutative by the definition of the homomorphism $\pi^{a,b}$ in (\ref{5.3}) and definitions of  differentials in (\ref{5.2}) and (\ref{2.1}).  
Let $w\in \Theta_n^{[a,b]}$,   then 
$w\in  \mathcal A_n^{[a,b]}(G)$ and $d(w) \in \mathcal A_{n-1}^{[a,b]}(G)$. By the above, we have 
$$
\pi^{[a,b]}(w) \in  \mathcal A_{n-2}(G) \ \ \text{and} \  \ 
\pi^{[a,b]}\left(d(w)\right)\in  \mathcal A_{n-3}(G).
$$
Furthermore,   by commutativity of  diagram (\ref{5.17}),  we have 
$$
\partial \left(\pi^{[a,b]} \left(w\right) \right)= \pi^{[a,b]}\left(d(w)\right)\in  \mathcal A_{n-3}(G).
$$
Hence 
$\pi^{[a,b]}(w)\in \Pi_{n-2}(G)$ 
and the statement of the theorem follows in the case (i).  The proof in the case (ii) is similar. 
\end{proof} 

\begin{corollary}\label{c5.9} Under assumptions of Theorem \ref{t5.8}  we have the following homomorphisms. In the case (i) ,  the morphism  of chain complexes in (\ref{5.12})   induces a homomorphism
$
\mathcal H^{[a,b]}_n(G) \to \mathcal  H_{n-2}(G) 
$
of primitive path  homology groups.  In the case (ii)  the morphism  of chain complexes in (\ref{5.13})   induces a homomorphism
$
\mathcal H^{[a,b]}_n(G) \to \widetilde{\mathcal  H}_{n-2}(G) 
$
of primitive path  homology groups.
\end{corollary} 

Now we  illustrate the relations of cluster homology  theory to geometric properties of digraphs. 

\begin{theorem}\label{t5.10} Let $G=(V,E)$ be  a digraph and $\overline G$ be the inverse directed digraph (see Definition \ref{d3.10}). For any two vertices  $v, w\in V$ and  $n\geq 0$, there is an isomorphism of chain complex 
$$
\Theta_*^{[a,b]}(G)\cong  \Theta_*^{[b,a]}(\overline G)
$$
and hence,  an isomorphism of primitive cluster homology groups
$$
\mathcal H_n^{[a,b]}(G)\cong \mathcal H_n^{[b,a]}(\overline G). 
$$
\end{theorem}
\begin{proof} Similar to  the proof of Theorem \ref{t3.11}. 
\end{proof}

\begin{definition}\label{d5.11} \rm Let  $G=(V,E)$ be a digraph. 
The \emph{directed suspension} $S^dG$  of the digraph $G$ is a digraph obtained from $G$ by adding two new vertices $a$ and $b$  and new arrows  $(a\to v), (v\to b)$  for all vertices $v\in V$. 
\end{definition}

\begin{theorem}\label{t5.12} Let $G=(V,E)$ be a digraph and  
$S^dG$ be the directed suspension defined by the vertices $a, b\notin V$. Then we have a natural isomorphism of 
 chain complexes 
$$
\Theta_*^{[a,b]}(S^dG)\cong \Pi_*(G)
$$
and, hence,  
$$
\mathcal H_n^{[a,b]}(S^dG)=\mathcal H_{n-2}(G). 
$$
\end{theorem}
\begin{proof}   For $n\geq 2$, we define a homomorphism 
$$
\pi\colon  \Lambda_n^{[a,b]}\left(V_{S^d(G)}\right)\to \Lambda_{n-2}(V)
\ \  \text{by}  \ \ \pi(e_{ai_{0}\dots i_{n-2}b})=e_{i_{0}\dots i_{n-2}}, 
$$  
 and we set $\pi=0\colon 0\to 0$ for $n=1$. It follows from Definition 
\ref{d5.11} that $\pi$ is an isomorphism and  its  restriction to 
$\mathcal A_{n}^{[a,b]}(S^dG)$ is an isomorphism 
$$
\mathcal A_{n}^{[a,b]}(S^dG)\to  \mathcal A_{n}(G).
$$
By definitions of  differentials $d$ in (\ref{5.2}) and $\partial$ in (\ref{2.1}) we obtain that the diagram  
\begin{equation}\label{5.16} 
\begin{matrix}
\Lambda_n^{[a,b]}\left(V_{S^dG}\right)& \overset{\pi}\longrightarrow & \Lambda_{n-2}(V)  \\
\\
d\downarrow&&\downarrow \partial  \\
\\
\Lambda_{n-1}^{[a,b]}\left(V_{S^dG}\right)& \overset{\pi}\longrightarrow 
&  \Lambda_{n-3}(V)\\
\end{matrix}
\end{equation} 
is commutative,   and  the result follows. 
\end{proof} 
\begin{example}\label{e5.13} \rm Consider a cyclic digraph $G=(V,E)$ below 
$$
\begin{matrix}
&& 2\ \\
& \swarrow &\uparrow\ \\
0& \to & 1.\\
\end{matrix}
$$
Let $S^dG$ be the directed suspension  defined by $a,b\notin V$ and 
$\mathbb Z$ be the ring of coefficients. 
By Corollary \ref{c3.8} and \cite[Example 2.7]{MiHomotopy}, 
$
\mathcal H_i(G)=H_i(G)=\mathbb Z
$ for $i=0,1$. Hence,   by Theorem \ref{t5.12}, 
$$
\mathcal H_3^{[a,b]}(S^dG)=\mathcal H_2^{[a,b]}(S^dG)=\mathbb Z.
$$
Now consider a digraph $H=(V_H, E_H)$ which 
has the same set of vertices as $S^dG$,  and the set of arrows $E_H$ is obtained from  $E_{S^dG}$ by deleting two  arrows $(a\to 2)$ and $2\to b$.  For $i\leq 4$, we have 
\begin{equation}\label{5.17}
\mathcal A_i^{[a,b]}(H)=\begin{cases}
0 & \ \ \text{for} \  \ i=0,1, \\
\langle e_{a0b},  e_{a1b}\rangle&   \ \ \text{for} \  \ i=2, \\
\langle e_{a01b}\rangle&   \ \ \text{for} \  \ i=3, \\
\langle  e_{a120b}\rangle&   \ \ \text{for} \  \ i=4.\\
\end{cases}
\end{equation}
The further calculation of the groups $\mathcal A_i^{[a,b]}(H)$ is similar. To find the groups $\Theta_*^{[a,b]}(H)$ we need to compute differentials 
$$
d\colon  \mathcal A_i^{[a,b]}(H)\to  \Lambda_i^{[a,b]}(H).
$$
It follows from  (\ref{5.17} that  $d=0$  for $i=1,2$  and  for $i=3,4$ we have:
\begin{equation}\label{5.18}
\begin{matrix} 
\text{for} \  i=3: & d(e_{a01b}) = e_{a1b}-e_{a0b},  \\
\text{for} \  i=4: & d(e_{a120b}) = e_{a20b}-e_{a10b}+ e_{a12b}.  \\
\end{matrix}
\end{equation}
It follows from (\ref{5.17}) and ({5.18}) that in dimensions $1\leq i\leq 4$ there is only one nontrivial homology group 
$\mathcal H_2^{[a,b]}(H)=\mathbb Z$. It is an easy exercise for readers to check similarly above  that  $\mathcal H_2^{[a,b]}(H)=0$ for $i\geq 5$. 
\end{example} 

\section{Primitive  homology   of paths   with a  fixed tail  or  a head vertex}\label{S6}
\setcounter{equation}{0}

In this section, for a digraph $G$, we construct primitive path homology theories  of paths   with a fixed tail  or a fixed   head vertex.   
The construction in both cases is similar to the previous section.
We describe relations between introduced primitive path  homology  and  the   path homology which was introduced in Sections \ref{S3} and \ref{S5}. 
We present several examples. In this section we fix a unitary and commutative ring $R$ of coefficients. 

\begin{definition}\label{d6.1}\rm  Let  $V$ be a finite set of vertices,   $a\in V$  be a vertex, and  $n\geq 0$.

 (i) A  path  in $\Lambda_n(V)$  is an $a$-\emph{tail path} if it is a linear combination of $n$-paths with the tail vertex $a$.  Such  path  is denoted $w^{[a,\cdot]}$. An elementary path fitting into the path $w^{[a,\cdot]}$ is called an \emph{elementary}  $a$-\emph{tail path}. 

(ii) A  path in $\Lambda_n(V)$  is an $a$-\emph{head path} if it is a linear combination of paths with the head vertex $a$.  Such  path  is denoted 
$w^{[\cdot,a]}$. An elementary path fitting into the path $w^{[\cdot,a]}$ is called an \emph{elementary}  $a$-\emph{head path}. 
\end{definition} 

This definition is a natural generalization of the definitions of paths in  a digraph $G$ which are given by 
   relations (\ref{4.1}) and  (\ref{4.2}). 

Let $n\geq 0$ and $a\in V$.  Denote by  $\Lambda _{n}^{[a,\cdot]}=\Lambda _{n}^{[a,\cdot]}( V)$ a free module generated by    elementary
$a$-tail $n$-paths.     The elements of $\Lambda_{n}^{[a,\cdot]} $ are called $a$-\emph{tail $n$-paths }  in $V$.   Denote by   $\Lambda _{n}^{[\cdot, a]}=\Lambda _{n}^{[\cdot,a]}( V)$ a free module generated by   elementary  
$a$-\emph{head $n$-paths}.     The elements of $\Lambda_{n}^{[\cdot, a]} $ are called $a$-\emph{head $n$-paths}  in $V$.  We set  
$\Lambda_{-1}^{[a,\cdot ]}=\Lambda_{-1}^{[\cdot,a]}=0.
$

For $n\geq 0$, define  \emph{differentials} 
\begin{equation}\label{6.1}
d^t\colon \Lambda _{n}^{[a,\cdot ]}\rightarrow \Lambda _{n-1}^{[a,\cdot ]} \ \ \text{and} \ \  d^h\colon \Lambda _{n}^{[\cdot, a ]}\rightarrow
 \Lambda _{n-1}^{[\cdot, a ]} 
\end{equation} 
 by setting on     basic elements 
\begin{equation}\label{6.2}
\begin{matrix}
 d^t(e_{ai_{0}\dots  i_{n-1}} )=\begin{cases} \sum\limits_{m=0}^{n-1}\left( -1\right)^{m}e_{ai_{0}\dots \widehat{i_{m}}\dots i_{n-1}}& \text{for} \ n\geq 1,\\
 0& \text{for} \  n=0, 
\end{cases}\\
\\
 d^h(e_{i_{0}\dots  i_{n-1}a} )=\begin{cases} \sum\limits_{m=0}^{n-1}\left( -1\right)^{m}e_{i_{0}\dots \widehat{i_{m}}\dots i_{n-1}a}& \text{for} \ n\geq 1,\\
 0& \text{for} \  n=0.
\end{cases}\\
\end{matrix}
\end{equation}
Let  $\widetilde{\Lambda}_*(V)$ be the augmented chain complex  defined in Section \ref{S2}. Recall that  $\widetilde{\Lambda}_n(V)=\Lambda_n(V)$  for 
$n\geq 0$,  $\widetilde{\Lambda}_{-1}(V)=R$,  and 
$\widetilde{\Lambda}_{-2}(V)=0$.

For $n\ge 0$, define  isomorphisms  of modules 
\begin{equation}\label{6.3}
\begin{matrix}
\pi^{[a,\cdot ]}\colon \Lambda_n^{[a,\cdot ]}(V)\to \widetilde \Lambda_{n-1}(V), \\
\\
\pi^{[\cdot, a ]}\colon \Lambda_n^{[\cdot, a]}(V)\to \widetilde \Lambda_{n-1}(V), \\
\end{matrix}
\end{equation}
by setting for $\pi^{[a,\cdot ]}$: 
$\pi^{[a,\cdot]}(e_{ai_{0}\dots i_{n-1}})=e_{i_{0}\dots i_{n-1}}$ for $n\geq 1$,  
$\pi^{[a,\cdot ]}(e_{a})= 1\in R$
 for $n= 0$, and $\pi^{[a,\cdot]}=0$  for $n=-1$.  For $\pi^{[\cdot, a ]}$ we set similarly: $\pi^{[\cdot,a]}(e_{i_{0}\dots i_{n-1}a})=e_{i_{0}\dots i_{n-1}}$ for $n\geq 1$,  
$\pi^{[\cdot,a]}(e_{a})= 1\in R$
 for $n= 0$, and $\pi^{[\cdot, a]}=0$  for $n=-1$.

\begin{proposition}\label{p6.2}  For $n\geq 0$, the following diagrams 
\begin{equation}\label{6.4} 
\begin{matrix}
\Lambda_n^{[a,\cdot ]}(V)& \overset{\pi^{[a,\cdot]}}\longrightarrow & \widetilde \Lambda_{n-1}(V)\\
\\
d^t\downarrow&&\downarrow \partial\\
\\
\Lambda_{n-1}^{[a, \cdot]}(V)& \overset{\pi^{[a,\cdot]}}\longrightarrow 
& \widetilde \Lambda_{n-2}(V) \\
\end{matrix}
\end{equation}
and
\begin{equation}\label{6.5} 
\begin{matrix}
\Lambda_n^{[\cdot,a ]}(V)& \overset{\pi^{[\cdot,a]}}\longrightarrow & \widetilde \Lambda_{n-1}(V)\\
\\
d^h\downarrow&&\downarrow \partial\\
\\
\Lambda_{n-1}^{[\cdot, a]}(V)& \overset{\pi^{[\cdot, a]}}\longrightarrow 
& \widetilde \Lambda_{n-2}(V)\\
\end{matrix}
\end{equation} 
are  commutative.
\end{proposition}
\begin{proof} Follows from the definition of  isomorphisms $\pi^{[a,\cdot]}$,  $\pi^{[\cdot,a]}$,    the  definition of  differential $\partial$ in (\ref{2.1}),  and definitions  of differentials $d^t, d^h$ in (\ref{6.2}). 
\end{proof} 

\begin{corollary}\label{c6.3} For $n\geq 0$, the modules $\Lambda_n^{[a, \cdot]}$ and  $\Lambda_n^{[\cdot, a]}$ with  differentials $d^t$ and $d^h$, respectively,  form chain complexes.
\end{corollary} 
\begin{proof} In digrams (\ref{6.4}) and (\ref{6.5}) the modules  
$\widetilde \Lambda_{*}(V)$ with the differential $\partial$ form  chain complexes.  The horizontal maps are isomorphisms. Now the statement follows from commutativity of these diagrams. 
\end{proof}

\begin{definition}\label{d6.4} \rm Let $G=(V, E)$ be a digraph,  $a\in V$ be a vertex,  and $n\geq 0$. 

(i) Let  $\mathcal A_n^{[a, \cdot ]}=\mathcal A_n^{[a, \cdot]}(G)$ be  a module generated 
by elementary allowed $a$-tail $n$-paths in $G$. We set $\mathcal A_{-1}^{[a, \cdot ]}=0$ and  $\mathcal A_n^{[a, \cdot]}=0$ if there are no such paths.  

(ii) Let  $\mathcal A_n^{[\cdot,a ]}=\mathcal A_n^{[\cdot, a]}(G)$ be  a module generated 
by elementary allowed $a$-head $n$-paths in $G$. We set $\mathcal A_{-1}^{[\cdot, a ]}=0$ and  $\mathcal A_n^{[\cdot,a]}=0$ if there are no such paths.  
\end{definition} 

It follows from Definition \ref{d6.4} that  for every digraph $G$ we have 
\begin{equation}\label{6.6}
\mathcal A_0^{[a, \cdot]}(G)=\mathcal A_0^{[\cdot, a]}(G)=\langle e_a\rangle  \ \text{for all} \  a \in V,
\end{equation} 
\begin{equation}\label{6.7}
\begin{matrix}
  \mathcal A_1^{[a, \cdot]}(G)=\begin{cases} \langle e_{ab}\, | \, (a\to b)\in E\rangle, &  \text{if arrows}\ (a\to b)\in E \ \text{exist}, \\
                                                                0,&  \text{otherwise},   \\
\end{cases}\\
  \mathcal A_1^{[\cdot, a]}(G)=\begin{cases} \langle e_{ba}\, |\,  (b\to a)\in E \rangle, &  \text{if arrows}\ (b\to a)\in E\  \text{exist}, \\
                                                                0,&  \text{otherwise.}  \\
\end{cases}\\
 \end{matrix}
\end{equation} 

For $n\geq -1$, the modules $\mathcal A_n^{[a, \cdot]}(G)$ and $\mathcal A_n^{[\cdot, a]}(G)$ are submodules of $\Lambda_{n}(V)$. For $n\geq 0$, we define submodules
\begin{equation}\label{6.8}
\begin{matrix}
\Theta_n^{[a,\cdot]}=\Theta_n^{[a,\cdot]}(G)=\{w\in \mathcal A_n^{[a, \cdot]}(G)\, |
\, d^t(w)\in  \mathcal A_{n-1}^{[a, \cdot]}(G)\}\\
\\
\Theta_n^{[\cdot, a]}=\Theta_n^{[\cdot, a]}(G)=\{w\in \mathcal A_n^{[\cdot, a]}(G)\, |
\, d^h(w)\in  \mathcal A_{n-1}^{[\cdot,a]}(G)\}\\
\end{matrix}
\end{equation} 
where $d^t$ and $d^h$ are differentials defined in 
(\ref{6.2}). We set $\Theta_{-1}^{[a,\cdot]}=\Theta_{-1}^{[\cdot,a]}=0$. 

\begin{proposition}\label{p6.5} The modules $\Theta_n^{[a,\cdot ]}(G) \, (n\geq 1) $   with the restriction of the differential $d^t$   and the modules 
$\Theta_n^{[\cdot, a]}(G), \, (n\geq -1)$ with the restriction of the differential $d^h$ form    chain complexes $\Theta_*^{[a,\cdot]}=\Theta_*^{[a,\cdot]}(G)$ and $\Theta_*^{[\cdot, a]}=\Theta_*^{[\cdot, a]}(G)$, respectively
\end{proposition} 
\begin{proof} Similar to Proposition \ref{p5.5}. 
\end{proof} 

 By  (\ref{6.8}) and  (\ref{6.7}),  and by the definition of differentials $d^h$ and $d^t$,  we obtain  that for any digraph $G=(V,E)$ 
\begin{equation}\label{6.9}
\Theta_i^{[a,\cdot]}(G)=\mathcal A_i^{[a,\cdot]}, \  \Theta_i^{[\cdot. a]}(G)=\mathcal A_i^{[\cdot, a]}  \ \ \text{for} \ \  -1\leq i\leq 1. 
\end{equation} 

\begin{definition}\label{d6.6}  \rm (i) The homology groups of the chain complex 
$\Theta_*^{[a,\cdot]}(G)$ are called \emph{primitive $a$-tail homology groups} of the digraph $G$ and are denoted by $\mathcal H_n^{[a,\cdot ]}(G)=\mathcal H_n^{[a,\cdot]}(G,R)\, (n\geq 0)$.

(ii)  The homology groups of the chain complex 
$\Theta_*^{[\cdot,a]}(G)$ are called \emph{primitive $a$-head homology groups} of the digraph $G$ and are denoted by $\mathcal H_n^{[\cdot, a]}(G)=\mathcal H_n^{[\cdot, a]}(G,R)\, (n\geq 0)$.
\end{definition} 

\begin{theorem}\label{t6.7} Let $G$ be a digraph and $\overline{G}$ be the digraph which is obtained by changing directions of all arrows. Then, for $n\geq 0$, there is an isomorphism $\mathcal H_n^{[\cdot, a]}(G)\cong \mathcal H_n^{[a, \cdot]}(G)$.
\end{theorem}
\begin{proof} Similar to the proof of Theorem \ref{t3.11}. 
\end{proof}

\begin{example}\label{e6.8} \rm Let $\mathbb Z$ be the ring of coefficients. 

(i) Consider a digraph $G$ below 
$$
\begin{matrix}
3&&2 &  &\\
\uparrow&\nearrow& &\searrow& \\
0& \to &1&\to &4.\\
\end{matrix}
$$
We have  
$
\mathcal A^{[0,\cdot]}_{-1}=0, \ \ \mathcal A^{[0,\cdot]}_0=\langle e_0\rangle, \  \ \mathcal A^{[0,\cdot]}_1=\langle e_{01}, e_{02}, e_{03}\rangle,
\ \ \mathcal A^{[0,\cdot ]}_2=\langle e_{014}, e_{024} \rangle,
$
and
$
\mathcal A^{[0,\cdot ]}_n=0  \ \text{for} \  n\geq 3.
$
In dimensions  one and two, the differential  $d^t$  is given as follows 
$$
d^t(e_{01})= d^t(e_{02})= d^t(e_{03})=e_0,
$$
$$
d^t(e_{014})= e_{04}-e_{01}, \ \ 
d^t(e_{024})= e_{04}-e_{02}. 
$$
Hence, 
$$
\Theta_0^{[0,\cdot]}(G)= \langle e_0\rangle,  \ \Theta_1^{[0,\cdot]}(G)= \langle e_{01}, e_{02}, e_{03}\rangle, \ \Theta_2^{[0,\cdot]}(G)= \langle  e_{014} - e_{024}\rangle,   
$$
and $\Theta_n^{[0,\cdot]}(G)=0$  for $n\geq 3$.
Now the direct computation gives  $
\mathcal H_1^{[0,\cdot]}(G)=\mathbb Z$ and  
$
\mathcal H_n^{[0,\cdot]}(G)=0  
$ for all $n\ne 1$.

(ii) Consider a digraph $G$ below 
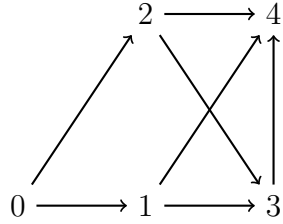
\begin{figure}[H]
\centering
\begin{tikzpicture}[scale=0.85]
\node (1) at (4,3) {$0$};
\node (2) at (6,3) {$1$};
\node (3) at (8,3) {$3$};
\node (4) at (6,6) {$2$};
\node (5) at (8,6) {$4$};

\draw (1) edge[ color=black!120, thick,-> ] (2);
\draw (1) edge[ color=black!120, thick, ->] (4);
\draw (2) edge[ color=black!120, thick, ->] (3);
\draw (2) edge[ color=black!120, thick, ->] (5);
\draw (4) edge[ color=black!120, thick,-> ] (5);
\draw (4) edge[ color=black!120, thick,->] (3);
\draw (3) edge[ color=black!120, thick,->] (5);
\end{tikzpicture}  

\end{figure}
We have  
$$
\begin{matrix}
\mathcal A^{[0,\cdot]}_{-1}=0, \ \ \mathcal A^{[0,\cdot]}_0=\langle e_0\rangle, \  \ \mathcal A^{[0,\cdot]}_1=\langle e_{01}, e_{02}\rangle,
\ \ \mathcal A^{[0,\cdot ]}_2=\langle e_{013}, e_{014}, e_{023}, e_{024} \rangle,
\\
\\
 \ \mathcal A^{[0,\cdot ]}_3=\langle e_{0134}, e_{0234}\rangle,  \ \ 
\mathcal A^{[0,\cdot ]}_n=0  \ \text{for} \  n\geq 4.\\
\end{matrix}
$$
In dimensions  $1\leq n\leq 3$, the differential  $d^t$  is given as follows 
$$
d^t(e_{01})= d^t(e_{02})=e_0,
$$
$$
d^t(e_{013})= e_{03}-e_{01}, \
d^t(e_{014})= e_{04}-e_{01}, \
d^t(e_{023})= e_{03}-e_{02}, \
d^t(e_{024})= e_{04}-e_{02}, 
$$
$$
\begin{matrix}
d^t(e_{0134})= e_{034}-e_{014}+e_{013}, \\
d^t(e_{0234})= e_{034}-e_{024}+e_{023}. \\
\end{matrix} 
$$
Hence, 
$$
\Theta_0^{[0,\cdot]}(G)= \langle e_0\rangle,  \ \Theta_1^{[0,\cdot]}(G)= \langle e_{01}, e_{02}\rangle, \ \Theta_2^{[0,\cdot]}(G)= \langle e_{013} - e_{023}, e_{014} - e_{024}\rangle, 
$$
$$
\Theta_3^{[0,\cdot]}(G)= \langle e_{0134} - e_{0234}\rangle,  \  \text{and}  \ \ \Theta_n^{[0,\cdot]}(G)=0 \ \ \text{for} \  \  n\geq 4.
$$
Now the direct computation gives  
$
\mathcal H_n^{[0,\cdot]}(G)=0  
$ for all $n$.

(iii) Consider a digraph  $G$  below 
$$
\begin{matrix}
 &&3& \to &4\\
&\nearrow&\uparrow&& \uparrow \\
1&\leftarrow&0&\to &2\\
\end{matrix}
$$
We have  
$$
\begin{matrix}
\mathcal A^{[0,\cdot]}_{-1}=0, \ \mathcal A^{[0,\cdot]}_0=\langle e_0\rangle\cong \mathbb Z, \ \mathcal A^{[0,\cdot]}_1=\langle e_{01}, e_{02},e_{03} \rangle \cong \mathbb Z^3,\\
\\
\mathcal A^{[0,\cdot ]}_2=\langle e_{013},  e_{024}, e_{034}  
\rangle \cong \mathbb Z^3, \ \mathcal A^{[0,\cdot ]}_3=\langle e_{0134}  
\rangle \cong \mathbb Z,\
\mathcal A^{[0,\cdot ]}_n=0  \ \text{for} \  n\geq 4.
\end{matrix}
$$
Now we can compute the differential  $d^t\colon  \mathcal A^{[0,\cdot ]}_n\to \Lambda_{n-1}^{0, \cdot]}$. We have 
$
d^t(e_0)=0$, 

$$
\begin{matrix}
d^t(e_{013})=e_{03}-e_{01},\\
d^t(e_{024})=e_{04}-e_{02},\\
d^t(e_{034})=e_{04}-e_{03},\\
\end{matrix}
$$
$$
d^t(e_{0134})=e_{034}-e_{014}+e_{013}. 
$$
Hence $\Theta_n^{[0,\cdot]}=0$ for $n\geq 3$, 
$\Theta_2^{[0,\cdot]}\langle    e_{013}, e_{024}-e_{034}\rangle \cong \mathbb Z^2$,  $\Theta_1^{[0,\cdot]}=\langle e_{01}, e_{02},e_{03} \rangle \cong \mathbb Z^3$, $\Theta_0^{[0,\cdot ]}=\langle e_0\rangle =\mathbb Z$, and $\Theta_n^{[0,\cdot ]}=0$ for $n=-1$. Thus we obtain 
$
\mathcal H_n^{[0,\cdot]}=0  
$ for all $n$. 
\end{example}

The notions of a cone  and of a suspension of a digraph were introduced in \cite{Mi2012},  \cite[Definition 4.20]{Mi2}. We use these digraph constructions to describe geometric properties of  primitive $a$-head  and $a$-tail homology. At first we recall these definitions. 

\begin{definition}\label{d6.9} \rm Let  $G=(V,E)$ be a digraph. 

(i) The \emph{cone} $CG$ 
 is obtained from the digraph $G$ by adding
a new vertex $a$   and new arrows  $(v\to a)$  for all vertices $v\in V$.

(ii) The \emph{inverse cone} $C^-G$ is obtained from the digraph $G$  by adding
a new vertex $a$   and new arrows  $(a\to v), (b\to v)$  for all vertices $v\in V$. 

(iii)  The \emph{suspension} $SG$ is obtained from the digraph $G$ by adding
two new vertices $a$ and $b$  and new arrows  $(v\to a), (v\to b)$  for all vertices $v\in V$. 

(iv) The \emph{inverse suspension} $S^-G$  is obtained from the digraph $G$ by adding
two new vertices $a$ and $b$  and new arrows  $(a\to v), (b\to v)$  for all vertices $v\in V$. 
\end{definition} 

\begin{theorem}\label{t6.10} Let $G=(V,E)$ be a digraph  and $CG$, $C^-G$ be  direct and  inverse cones which are  defined by the vertex $a\notin V$. Then there are the following relations between primitive path homology groups: 
$$
\mathcal H_n^{[v,a]}(CG)= \mathcal H_{n-1}^{[v, \cdot]}(G) \ \ \text{for every} \ \   v\in V, 
$$
$$
\mathcal H_n^{[a, v]}(C^-G)= \mathcal H_{n-1}^{[\cdot, v]}(G) \ \ \text{for every} \ \   v\in V. 
$$
Let $SG$, $S^-G$ be direct and  inverse suspensions  which are  defined by the vertices  $a, b\notin V$. Then there are the following relations between primitive path homology groups: 
$$
\mathcal H_n^{[v, a]}(SG)= \mathcal H_n^{[v, b]}(SG)=\mathcal H_n^{[v,a]}(CG)=\mathcal H_{n-1}^{[v, \cdot]}(G)  \ \ \text{for every} \ \ v\in V, 
$$
$$
\mathcal H_n^{[a, v]}(S^-G)= \mathcal H_n^{[b, v]}(S^-G)= \mathcal H_n^{[a,v]}(C^-G)=\mathcal H_{n-1}^{[ \cdot, v]}(G)  \ \ \text{for every} \ \ v\in V. 
$$
\end{theorem}
\begin{proof}  Similar to the proof of Theorem \ref{t5.12}. 
\end{proof} 

\section{Functoriality}\label{S7}
\setcounter{equation}{0}

In this section, we prove functoriality of primitive  homology  theories constructed in Sections  \ref{S5} and \ref{S6} on  the category $\mathcal{AD}$ of asymmetric digraphs.

\begin{theorem}\label{t7.1}  Let 
$f\colon G\to H$ be a  map of  asymmetric digraphs. For two vertices $a,b\in V_G$  the homomorphism  
$f_{\sharp}\colon \Lambda_n(V_G)\to \Lambda_n(V_H) \, (n\geq -1)$ defined in (\ref{3.4}) 
provides morphisms of chain complexes:
\begin{equation}\label{7.1}
\begin{matrix}
f_{*}\colon \Theta_*^{[a,b]}(G)\to \Theta_*^{[f(a),f(b)]}(H),\\
\\
f_{*}\colon \Theta_*^{[a,\cdot]}(G)\to \Theta_*^{[f(a),\cdot]}(H),\\
\\
f_{*}\colon \Theta_*^{[\cdot,b]}(G)\to \Theta_*^{[\cdot,f(b)]}(H)
\end{matrix}
\end{equation}
which  induce homomorphisms of primitive homology groups:
\begin{equation}\label{7.2}
\begin{matrix}
f_{*}\colon \mathcal H_*^{[a,b]}(G)\to  \mathcal H_*^{[f(a),f(b)]}(H),\\
\\
f_{*}\colon \mathcal H_*^{[a,\cdot]}(G)\to  \mathcal H_*^{[f(a),\cdot]}(H),\\
\\
f_{*}\colon \mathcal H_*^{[\cdot,b]}(G)\to  \mathcal H_*^{[\cdot,f(b)]}(H). 
\end{matrix}
\end{equation}
\end{theorem} 
\begin{proof} Similar to the proof of Proposition \ref{p3.5}. 
\end{proof}

\bibliographystyle{amsplain}
\bibliography{biblio}

Jingyan Li: \emph{Beijing 
Institute of Mathematical
Sciences and Application,  No. 544 Hefangkou Village, Huairou, Beijing 
 China.}

e-mail: jingyanli@bimsa.cn

Yuri Muranov: \emph{Faculty of Mathematics and Computer Science, University
of Warmia and Mazury in Olsztyn, Sloneczna 54 Street, 10-710 Olsztyn, Poland. }

e-mail: muranov@matman.uwm.edu.pl

Jie Wu:\emph{Beijing 
Institute of Mathematical
Sciences and Application,  No. 544 Hefangkou Village, Huairou, Beijing 
 China.}

e-mail: wujie@bimsa.cn

Shing-Tung Yau:  \emph{Room 215, Yau Mathematical Sciences Center, Jing Zhai, Tsinghua University, Hai District, Beijing,  China.}

e-mail: styau@tsinghua.edu.cn

\end{document}